\documentclass[reqno,10pt]{amsart}
\usepackage{amscd,amssymb,amsmath,color,stmaryrd}
\usepackage{latexsym,wasysym}
\usepackage{accents}
\usepackage{mathrsfs}
\usepackage[all]{xy}
\usepackage{imakeidx}
\usepackage{braket}

\theoremstyle{definition}
\newtheorem{theor}[subsubsection]{Theorem}
\newtheorem{prop}[subsubsection]{Proposition}
\newtheorem{lem}[subsubsection]{Lemma}
\newtheorem{cor}[subsubsection]{Corollary}
\newtheorem{defin}[subsubsection]{Definition}
\newtheorem{rem}[subsubsection]{Remark}
\newtheorem{rems}[subsubsection]{Remarks}
\newtheorem{exam}[subsubsection]{Example}
\newtheorem{exams}[subsubsection]{Examples}

\def\hatotimes{{\widehat\otimes}}
\def\toisom{\widetilde{\to}}

\def\Hom{{\rm Hom}}
\def\End{{\rm End}}
\def\Ker{{\rm Ker}}

\def\Res{{\rm Res}}

\def\bfZ{{\bf Z}}
\def\bfQ{{\bf Q}}
\def\bfR{{\bf R}}
\def\bfa{{\bf a}}

\def\calF{{\mathcal F}}
\def\calW{{\mathcal W}}

\def\oV{{\overline V}}
\def\oj{{\overline j}}

\font\goth=eufm10
\def\gothg{{\hbox{\goth g}}}


\begin{document}

\title{Non-Archimedean Vertex Algebras}

\author{Victor G. Kac}

\maketitle

\setcounter{tocdepth}{2}

\begin{abstract}
Foundations of the theory of vertex algebras are extended to the non-Archimedean setting.
\end{abstract}

\section{Introduction}
This is a joint work with Vladimir G. Berkovich.

The notion of a vertex algebra was introduced by Borcherds [Bo] in terms of infinitely many bilinear products $a_{(n)}b$, $n\in\bfZ$, satisfying a cubic identity, similar to, but much more complicated than, the Jacobi identity for a Lie algebra.  This notion is a rigorous definition of the chiral part of a 2-dimensional conformal field theory, studied intensively by physicists and mathematicians since the landmark paper of Belavin, Polyakov and Zamolodchikov [BPZ].

Much simpler axioms were used in the book [K], where it was demonstrated that they are equivalent to the Borcherds identity.  These axioms are those of the paper [FKRW], which were inspired by the paper of Goddard [G].  They are formulated in terms of fields $a(z) = \sum_{n\in\bfZ}a_{(n)}z^{-n-1}$, which are the generating series of the operators of $n^{\mathrm{th}}$ multiplication $a_{(n)}$, where $a_{(n)}b=0$ for $n \gg 0$.

In more detail, a vertex algebra is a vector space $V$ (space of states) over a unital commutative ring $K$, endowed with a vector $\ket{0}$ (vacuum vector), an endomorphism $T$ (infinitesimal translation operator), and a linear map of $V$ to the space of fields (the state-field correspondence) $a\mapsto a(z) = \sum_{n\in\bfZ}a_{(n)}z^{-n-1}$, where $a_{(n)}\in \End(V)$, $a_{(n)}b=0$ for each $b\in V$, if $n \gg 0$, subject to the following axioms ($a$, $b\in V$):
\begin{itemize}
	\item [(A.1)]
({\it vacuum}) $T\ket{0}=0$, $\ket{0}(z)=I_V$, $a(z)\ket{0}\in a+zV[[z]]$;
    \item[(A.2)]
({\it translation covariance}) $[T, a(z)] = \frac{d}{dz}a(z)$;
    \item[(A.3)]
({\it locality}) $(z-w)^N [a(z), b(w)] = 0$ for $N\gg 0$.
\end{itemize}
In fact, as explained in [K], these axioms form the essential part of the Wightman axioms of an arbitrary QFT [W].

In the present paper we develop the notion of a non-Archimedean vertex algebra, along the lines of [K], which generalizes the notion of a vertex algebra.  See Definition \ref{def-vertex}, where they are called vertex $K$-algebras.  Namely, we take for $K$ a commutative non-Archimedean Banach ring and for $V$ a Banach $K$-module.  (In the case of ``ordinary'' vertex algebras one takes the trivial Banach norm on $K$ and on $V$.)  Furthermore the condition on $a(z)$ to be a field is replaced by
$$
\lim_{n \to +\infty} a_{(n)}b=0~\mathrm{for~each}~b \in V,
$$
the vacuum and translation covariance axioms remain the same, but the locality axiom generalizes to
$$
\lim_{N \to +\infty} (z-w)^N [a(z), b(w)] = 0.
$$
Finally, the field-state correspondence
$$
a(z) \mapsto a(z)\ket{0} |_{z=0}
$$
is not surjective in general, but it is injective, with image $V^\prime$ (Lemma \ref{lem-fs}), so that in general the state-field correspondence is an unbounded operator with the domain $V^\prime$.  Thus, we cannot identify the spaces of states and of fields (like in [FR]), as is always done in the case of ordinary vertex algebras.

A non-Archimedean vertex algebra $V$ is called admissible if the field-state homomorphism of Banach $K$-modules is admissible.  In this case $V^\prime$ is a closed submodule of $V$ and the state-field homomorphism is bounded.  In particular, if $K$ is a field with a non-trivial valuation, $V$ is admissible if and only if $V = V^\prime$.

In the first part of the paper we develop the theory of formal distributions and fields in the non-Archimedean Banach setting.  In particular, we prove the Decomposition Theorem for any local formal distribution in two indeterminates (Lemma \ref{lem-repdelta}); we define the $n^\mathrm{th}$ product $a(z)_{(n)}b(z)$ of two fields for each $n\in\bfZ$, and establish their properties, in particular, Dong's Lemma \ref{lem-Dong}.

In the second part we define a non-Archimedean vertex algebra over a non-Archimedean commutative Banach ring $K$, which is torsion free as an abelian group (Definition \ref{def-vertex}), and prove their basic properties.  In particular, we prove Goddard's uniqueness theorem (Lemma \ref{lem-Goddard}), the $n^\mathrm{th}$ product identity (Lemma \ref{prop-Lconf} ii), and the Borcherds identity (Corollary \ref{cor-Borcherds}).  Furthermore, we establish the Extension Theorem (Theorem \ref{theor-gener}), which along with the ring extension lemma (Lemma \ref{lem-fsurj}) allows us to construct examples of non-Archimedean vertex algebras.  In \S \ref{subsec-lieassoc} we discuss their construction, associated to Lie (super)algebras.  Based on this, we construct in Subsection \ref{subsec-bosfer} the free boson and fermion vertex algebras, and the universal Virasoro and affine vertex algebras.

Finally in \S \ref{subsec-lieconformal} we introduce the notion of a non-Archimedean Lie conformal algebra of radius $r>0$ (Definition \ref{defn-lieconformal}), and show that for a non-Archimedean vertex algebra the space of fields carries a structure of a non-Archimedean Lie conformal algebra of radius $r$, defined by the $\lambda$-bracket
$$
\left[a(z)_\lambda b(z)\right] = \sum_{n=0}^{\infty}\frac{\lambda^n}{n!}\left(a(z)_{(n)}b(z)\right),
$$
where
$$
\lim_{n\to +\infty}r^n ||a(z)_{(n)}b(z)|| / n! = 0,
$$
see Proposition \ref{prop-Lconf}.  Here we assume that the base Banach ring $K$ contains $\bfQ$, and we let $r=1$ (resp. $r=|p|^\frac{1}{p-1}$) if the valuation $|.|$ of $\bfQ$, induced from $K$ is trivial (resp. $p$-adic).

When our paper was nearing completion, we learned about the paper [FM], where $p$-adic vertex algebras are introduced and studied.  Their context is more restrictive, since they assume that the base ring $K$ is a $p$-adic field and that the state-field correspondence is bijective.  On the other hand, they consider a connection of their theory to $p$-adic modular forms.
	
\section{Formal distributions and fields}

\subsection{Non-Archimedean Banach algebras and modules (see [B] for details)}

A {\it non-Archimedean Banach ring} is a unital ring $A$ complete
with respect to a non-Archimedean Banach norm. Recall that the
latter is a function $||\ ||:A\to \bfR_+$ satisfying the following
axioms: (1) $||a||=0$ if and only if $a=0$; (2) $||a b||\le
||a||\cdot||b||$; and (3) $||a+b||\le\max\{||a||,||b||\}$. An
example of a non-Archimedean norm is the {\it trivial} one $|\
|_0$ which takes value $1$ at all nonzero elements. Any abstract
ring is complete with respect to the trivial norm, i.e., it can be
considered as a non-Archimedean ring.

A {\it Banach $A$-module} is a left $A$-module $V$ complete with
respect to a {\it Banach norm}, i.e., a function $||\
||:V\to\bfR_+$ satisfying the following axioms: (1) $||v||=0$ if
and only if $v=0$; (2) $||a v||\le||a||\cdot ||v||$ for all $a\in
A$; and (3) $||u+v||\le\max\{||u||,||v||\}$. A {\it Banach
	$A$-algebra} is a non-Archimedean Banach ring $B$ provided with a
homomorphism of rings $\varphi:A\to B$ such that $||\varphi(a)||
\le ||a||$ for all $a\in A$.

{\it Throughout the text, K is a fixed commutative non-Archimedean Banach ring, and all modules considered are K-modules or Banach K-modules, unless otherwise stated.}

Given Banach $K$-modules $U$ and $V$, a $K$-linear operator
$T:U\to V$ is said to be {\it bounded} if there exists a positive
constant $C$ such that $||T u||\le C||u||$ for all $u\in U$. The
infinum of such $C$'s is called the {\it norm} of $T$ and denoted
by $||T||$. Banach $K$-modules form an additive category in which sets of morphisms are the spaces $\Hom(U,V)=\Hom_K(U,V)$ of all bounded
$K$-linear operators from $U$ to $V$, called {\it homomorphisms of Banach $K$-modules}. Such a space, provided with the above norm, is
a Banach $K$-module. If $V=K$, it is denoted by $U^\ast$ and, if
$U=V$, it is denoted by $\End(V)$. Notice that $\End(V)$ is a
Banach $K$-algebra, and to provide a Banach $K$-module $V$ with
the structure of a Banach $A$-module for a Banach $K$-algebra $A$
is the same as to provide $\End(V)$ with the structure of a Banach
$A$-algebra. Given Banach $A$-modules $U$ and $V$, the space
$\Hom_A(U,V)$ of bounded homomorphisms of $A$-modules is a closed
$K$-submodule of $\Hom(U,V)$. Of course, if $A$ is commutative, it
is again a Banach $A$-module.

Any closed $K$-submodule $U$ of a Banach $K$-module $V$ is a Banach $K$-module with respect to the induced norm. In this case, the quotient $V/U$ is a Banach $K$-module with respect to the quotient norm (the norm of a coset is the infimum of norms of elements from this coset). A homomorphism of Banach $K$-modules $\varphi:U\to V$ is said to be {\it admissible} if the Banach norm on ${\rm Im}(\varphi)$, induced from $V$, is equivalent to the norm induced from ${\rm Coim}(\varphi) = U/\Ker(\varphi)$ by the algebraic isomorphism ${\rm Coim}(\varphi) \toisom {\rm Im}(\varphi)$. In this case, ${\rm Im}(\varphi)$ is complete with respect to the norm, induced from $V$, i.e., it is closed in $V$ and the above algebraic isomorphism is in fact an isomorphism of Banach $K$-modules.  A bijective homomorphism of Banach $K$-modules is an isomorphism (in the category of Banach $K$-modules) if and only if it is admissible. Obviously, any isometry is an admissible homomorphism.

Given Banach $K$-modules $U$, $V$ and $W$, a $K$-bilinear map
$\varphi:U\times V\to W$ is {\it bounded} if there is a positive
constant $C$ with $||\varphi(u,v)|| \le C||u||\cdot||v||$ for all
$u\in U$ and $v\in V$. The {\it completed tensor product} of
Banach $K$-modules $U$ and $V$ is a Banach $K$-module $U\hatotimes_K
V$ provided with a bounded $K$-bilinear map $U\times V\to
U\hatotimes_K V$ such that any bounded $K$-bilinear map $U\times V
\to W$ goes through a unique bounded $K$-linear operator
$U\hatotimes_K V \to W$. It is clear that the completed tensor
product is unique up to a unique isomorphism, and it is
constructed as follows. Consider the usual tensor product
$U\otimes_K V$ (over $K$), and provide it with a real valued
function $||\ ||$ as follows: $||x|| = \inf \max\limits_i
||u_i||\cdot||v_i||$, where the infimum is taken over all
presentation of $x\in U\otimes V$ in the form of a finite sum
$\sum_i u_i\otimes v_i$. This function is a semi-norm (i.e., it
possesses the properties (2) and (3) of norms), and it gives rise
to a norm on the quotient of $U\otimes_K V$ by its kernel. Then
$U\hatotimes_K V$ is the completion of the quotient with respect to
that norm. Notice that, if $V$ is a Banach $K'$-module for a commutative non-Archimedean Banach $K$-algebra $K'$, then $U\hatotimes_K V$ has the structure of a Banach $K'$-module (defined in the evident way).

\begin{rem}\label{rem-val}
(i) The function $||\ ||':K\to\bfR_+$, defined by $||a||' = {\rm sup}_{b\not=0} \frac{||ab||}{||b||}$ (it is the norm of the operator of multiplication by $a$), is a Banach norm, equivalent to the norm $||\ ||$. For example, $||1||'=1$.  Similarly, for any Banach $K$-module $V$, the function $||\ ||':V\to\bfR_+$, defined by $||v||'={\rm sup}_{a\not=0} \frac{||av||}{||a||}$ is a Banach norm on $V$, equivalent to its norm $||\ ||$, and one has $||av||'\le||a||'\cdot||v||'$. Thus, we may always assume that the norm of elements from the image of $\bfZ$ in $K$ is at most one.

(ii) A non-Archimedean Banach norm on a commutative ring $A$ is said to be a (real) {\it valuation} if it is multiplicative, i.e., $||ab||=||a||\cdot||b||$ for all $a,b\in A$.
For example, by the Ostrowski theorem, any valuation on $\bfZ$ is either trivial, or $p$-adic for some prime integer $p$. The later means that $|p|<1$. In this case, $|n|=|p|^e$, where $e$ is the maximal power of $p$ such that $p^e$ divides $n$, and the completion of $\bfZ$ is the ring of $p$-adic integers $\bfZ_p$.
\end{rem}

\subsection{Examples of Banach modules and algebras}\label{subsec-Banmod}
For a Banach $K$-module $V$, let $V[[z]]$ denote the Banach $K$-module of all
formal power series $f=\sum_{n=0}^\infty v_n z^n$ with
$||f||=\sup\limits_n ||v_n||<\infty$, and let $V\{z\}$ denote the
Banach $K$-submodule of $V[[z]]$ consisting of the series with
$v_n\to0$ as $n\to\infty$. Similarly, let $V[[z^{\pm1}]]$ denote
the Banach $K$-module of all formal power series
$f=\sum_{n\in\bfZ} v_n z^n$ with $||f||=\sup\limits_n
||v_n||<\infty$, and let $V\{z^{\pm1}\}$ denote the Banach
$K$-submodule of $V[[z]]$ consisting of the series with $v_n\to0$
as $n\to\pm\infty$. We also introduce an intermediate Banach
$K$-submodule $V[[z]] \subset V((z)) \subset V[[z^{\pm1}]]$ that
consists of the series $f=\sum_{n\in\bfZ} v_n z^n \in
V[[z^{\pm1}]]$ with $v_n\to0$ as $n\to-\infty$.

Notice that for any Banach $K$-module $V$ the canonical bounded
$K$-linear operator $V\hatotimes K\{z\} \to V\{z\}$ is an
isomorphism, but the similar operators $V\hatotimes K[[z]] \to
V[[z]]$ and $V\hatotimes K((z)) \to V((z))$ are not.

One can iterate the above constructions. For example, both Banach
$K$-modules $V[[z^{\pm1}]][[w^{\pm1}]]$ and
$V[[w^{\pm1}]][[z^{\pm1}]]$ (resp. $V\{z^{\pm1}\}\{w^{\pm1}\}$ and
$V\{w^{\pm1}\}\{z^{\pm1}\}$) are canonically identified and
denoted by $V[[z^{\pm1},w^{\pm1}]]$ (resp.
$V\{z^{\pm1},w^{\pm1}\}$). But $V((z))((w))$ and $V((w))((z))$ are
different Banach $K$-submodules of $V[[z^{\pm1},w^{\pm1}]]$. Their
intersection, denoted by $V((z,w))$, consists of all series
$f=\sum_{m,n\in\bfZ} v_{m n} z^m w^n$ from $V[[z^{\pm1},w^{\pm1}]]$
with $v_{m n}\to0$ as $\min(m,n) \to-\infty$.

Assume now that the above $V$ is a Banach $K$-algebra $A$. Then
all of the introduced spaces are Banach $A$-modules. Moreover, the
spaces $A[[z]]$, $A\{z\}$, $A\{z^{\pm1}\}$ and $A((z))$ are Banach
$K$-algebras, $A((z))$ is a Banach $A\{z^{\pm1}\}$-algebra and, if
$V$ is a Banach $A$-module, $V[[z^{\pm1}]]$ and $V((z))$ are
Banach $A\{z^{\pm1}\}$-modules. Notice that although
$A[[z^{\pm1},w^{\pm1}]]$ is not an algebra, the product $f(z)
g(w)$ of any pair of elements $f(z)\in A[[z^{\pm1}]]$ and $g(w)\in
A[[w^{\pm1}]]$ is well defined, and one has $||f(z)g(w)||\le
||f(z)||\cdot||g(w)||$.

For example, if the norms on $K$ and $A$ are trivial,
$A\{z^{\pm1}\}$ is the algebra of Laurent polynomials $A[z^{\pm
	1}]$, $A((z))$ is the algebra of Laurent power series (i.e., those
with at most a finite number of nonzero coefficients at the
monomials of negative degree), and $A[[z^{\pm1}]]$ is the space of
bilateral formal power series in $z$.

\subsection{A Banach $K$-module $V((z,w))\{(z-w)^{-1}\}$}
Let $V$ be a Banach $K$-module.

\begin{lem}\label{lem-zminw}
The following is true:
\begin{itemize}
	\item [(i)]
the operator of multiplication by $z-w$ on $V[[w^{\pm1}]]((z))$ is an isometry;
    \item[(ii)]
$\cap_{n=0}^\infty (z-w)^n V((z,w))=0$;
    \item[(iii)]
the operator of multiplication by $u(z-w)-1$ on $V((z,w))\{u\}$ is an isometry.
\end{itemize}
\end{lem}

\begin{proof}
(i) If $f=\sum_{n\in\bfZ} f_n(w) z^{-n}\in
V[[w^{\pm1}]]((z))$, then $||f||=\sup\limits_n ||f_n||$ and
$f_n\to0$ as $n\to+\infty$. One has $(z-w)f = \sum_{n\in\bfZ}
(f_{n+1}-wf_n)z^{-n}$ and, therefore,
$$
||(z-w)f|| = \sup\limits_n ||f_{n+1}-wf_n||\ .
$$
Assume that $f\not=0$. Given $\varepsilon>0$, let $m$ be maximal
with the property $||f_m||>||f||-\varepsilon$ and $||f_n||\le
||f||-\varepsilon$ for all $n\ge m+1$. Then $||w f_m||=||f_m||$
and, therefore, $||f_{m+1} - w f_m|| = ||f_m||>||f||-\varepsilon$,
i.e., $||(z-w)f||>||f||-\varepsilon$. Since the latter is true for
every $\varepsilon>0$, the required fact follows.

(ii) For $n\in\bfZ$, let $L_n$ be the bounded $K$-linear operator
$$
V((z,w))\to V((z)): f=\sum_{i,j\in\bfZ} v_{ij} z^i w^j \mapsto
\sum_{i\in\bfZ} v_{i,n-i}z^i\ .
$$
One can easily see that $||f||= \sup\limits_n ||L_n f||$ and
$L_{n+1}((z-w)f)=(z-1)L_n(f)$. Thus, to verify the required fact,
it suffices to show that $\cap_{n=0}^\infty (z-1)^n V((z))=0$.
Given an element $f=\sum_{i\in\bfZ} v_i z^i\in V((z))$, let $l(f)$
and $r(f)$ be the minimal and maximal $i$ with $||v_i||=||f||$,
and set $c(f)=r(f)-l(f)$. One can easily see that $l((z-1)f)=l(f)$
and $r((z-1)f)=r(f)+1$ and, therefore, $c((z-1)f)=c(f)+1$. It
follows that an element $f\in V((z))$ can be divisible by at most
$c(f)$-th power of $z-1$.

(iii) If $f=\sum_{n=0}^\infty f_n(z,w) u^n\in V((z,w))\{u\}$, then
$f_n\to0$ as $n\to+\infty$ and $||f||=\max\limits_n ||f_n||$. One
has $(u(z-w)-1) f = -f_0 + \sum_{n=1}^\infty ((z-w)f_{n-1} - f_n)
u^n$ and, therefore
$$
||(u(z-w)-1) f|| = \max\{ ||f_0||,\max\limits_{n\ge1}
||(z-w)f_{n-1} -  f_n||\}\ .
$$
If $||f||=||f_0||$, the statement is trivial. Assume that
$||f||>||f_0||$, and let $m$ be maximal with the property
$||f||=||f_m||>||f_n||$ for all $n\ge m+1$. By (i), $||(z-w)f_m||
= ||f_m||$ and, therefore, $||(z-w) f_m - f_{m+1}||=||f_m||=||f||$. This implies the required fact.
\end{proof}

Let $V((z,w))\{(z-w)^{-1}\}$ denote the quotient of
$V((z,w))\{u\}$ by the closed $K$-submodule $(u(z-w)-1)V((z,w))\{u\}$.
It is a Banach $K$-module with respect to the quotient norm. We
denote by $(z-w)^{-1}$ the image of $u$ in it. Of course, if the
norm on $V$ is trivial, then $V((z,w))\{(z-w)^{-1}\} =
V((z,w))[(z-w)^{-1}]$.

\begin{cor}\label{cor-zminw}
Every element $f\in
V((z,w))\{(z-w)^{-1}\}$ has a unique representation in the form
$$
f(z,w)=h(z,w)+\sum_{n=0}^\infty {{g_n(w)}\over {(z-w)^{n+1}}}\ ,
$$
where $h\in V((u,w))$ and $g_n\in V((w))$ are such that $g_n\to0$
as $n\to\infty$. One also has $||f||=\max\{||h||,\max\limits_n
||g_n||\}$.
 \hfill$\square$	
\end{cor}

There are the following bounded homomorphisms of Banach
$K$-modules
$$
V((w))((z)) \buildrel {i_{z,w}} \over \longleftarrow
V((z,w))\{(z-w)^{-1}\} \buildrel {i_{w,z}} \over \longrightarrow
V((z))((w))\ ,
$$
whose restrictions to $V((z,w))$ are the canonical ones and which
take $(z-w)^n$ for $n\in\bfZ$ to its decompositions in the domains $|z|>|w|$ and $|z|<|w|$, respectively. For example, one has
$$
i_{z,w}\left({1\over{z-w}}\right) = \sum_{n=0}^\infty w^n z^{-n-1}{\ \
	\rm and\ \ }i_{w,z}\left({1\over{z-w}}\right) = -\sum_{n=-1}^{-\infty} w^n
z^{-n-1}\ .
$$
We denote by $\Delta$ the operator (see \S \ref{subsec-delta})
$$
V((z,w))\{(z-w)^{-1}\} \to V[[z^{\pm1},w^{\pm1}]]: f\mapsto
i_{z,w}(f)-i_{w,z}(f)\ .
$$
Notice that this operator commutes with multiplication by $z-w$.

\begin{rems}\label{rem-minusw}
(i) The correspondence $z\mapsto z$, $w\mapsto -w$ extends in the evident way to an isometric automorphism of the Banach $K$-module $V((z,w))$. This automorphism takes $z-w$ to $z+w$ and extends to an isometric isomorphism $V((z,w))\{(z-w)^{-1}\} \toisom V((z,w))\{(z+w)^{-1}\}: f(z,w) \mapsto f(z,-w)$. One can therefore apply the constructions of this subsection and, in particular, the operators $i_{z,w}$ and $i_{w,z}$ to the Banach $K$-module $V((z,w))\{(z+w)^{-1}\}$.  

(ii) The elements $i_{z,w}(f)$ and $i_{w,z}(f)$ are often denoted by $f_>$ and $f_<$, respectively.
\end{rems}

\subsection{Derivation and residue operators}\label{subsec-deres}

For a Banach $K$-module $V$, the derivation $\partial_z:V[[z^{\pm1}]] \to
V[[z^{\pm1}]]$ is a bounded $K$-linear operator of norm one. More
generally, for every $i\ge0$, there is a well defined bounded
$K$-linear operator of norm one $\partial_z^{(i)} = {1\over{i!}} \partial_z^i:
V[[z^{\pm1}]] \to V[[z^{\pm1}]]$
$$
\partial_z^{(i)}(\sum_{n\in\bfZ} v_n z^n) =
\sum_{n\in\bfZ} \binom{n}{i} v_n z^{n-i}\ .
$$
Notice that the operators $\partial_z^{(i)}$ commute
between themselves and preserve the Banach $K$-submodules
$V[[z]]$, $V((z))$ and $V\{z^{\pm1}\}$. In particular, they induce
operators on $V((z,w))$, and the latter are extended to bounded
$K$-linear operators on $V((z,w))\{(z-w)^{-1}\}$, which commute
with the operators $i_{z,w}$, $i_{w,z}$, and $\Delta$. Since ${1\over{(z-w)^{i+1}}} = \partial^{(i)}_\varpi ({1\over{z-w}})$, it follows that
 $$
i_{z,w}\left({1\over{(z-w)^{i+1}}}\right) = \sum_{n\ge0} \binom{n}{i} w^{n-i}z^{-n-1}\ {\rm and\ }  
 $$
 $$
i_{w,z}\left({1\over{(z-w)^{i+1}}}\right) = -\sum_{n\le-1} \binom{n}{i} w^{n-i}z^{-n-1}\ .
 $$

Furthermore, let $\Res_z$ denote the bounded $K$-linear operator
of norm one (the {\it residue} operator)
$$
V[[z^{\pm1}]]\to V: f=\sum_{n\in\bfZ} v_n z^n \mapsto
\Res_{z} f(z)  = v_{-1}\ .
$$
Notice that $\Res_z\circ \partial_z^{(i)}=0$ for all
$i\ge1$. If $V$ is a Banach $A$-module, the residue operator
defines a bounded $K$-bilinear pairing
$$
A\{z^{\pm1}\} \times V[[z^{\pm1}]]\to V: (\varphi,f) \mapsto
\Res_z(\varphi f)\ .
$$

\begin{lem}\label{lem-localdist}
The above pairing gives rise to an isometric isomorphism of Banach $K$-modules
$$
V[[z^{\pm1}]] \toisom \Hom_A(A\{z^{\pm1}\},V)\ .
$$
\end{lem}

\begin{proof} 
This homomorphism takes an element $v(z)=\sum_{j\in\bfZ} v_j z_j \in V[[z^{\pm1}]]$ to the homomorphism $\chi_v: A\{z^{\pm1}\}\to V$, which sends an element $\varphi(z)=\sum_{i\in\bfZ} a_i z^i\in A\{z^{\pm1}\}$ to the element $\sum_{i+j=-1} a_i v_j$. It follows that $||\chi_v(\varphi(z))||\le ||v(z)||\cdot ||\varphi(z)||$, i.e., $||\chi_v||\le||v(z)||$. The homomorphism in the opposite direction takes a homomorphism $\chi:A\{z^{\pm1}\}\to V$ to the element $v_\chi(z) = \sum_{j\in J} \chi(z^{-j-1}) z^j$. Both homomorphisms are inverse to each other, and $||v_\chi(z)|=\sup_j\{||\chi(z^{-j-1})||\} = ||\chi||$.
\end{proof}

Because of Lemma \ref{lem-localdist}, elements of $V[[z^{\pm1}]]$ are called {\it formal distributions}
on the space of {\it test functions} $A\{z^{\pm1}\}$. The operator
$\Res_z:V((z,w))\to V((w))$ is naturally extended to a bounded
homomorphism of $K((w))$-modules
$$
\Res_z: V((z,w))\{(z-w)^{-1}\} \to V((w)): f\mapsto \Res_z(i_{z,w}(f))\ .
$$
For example, $\Res_z({1\over{z-w}}) = 1$.

For $m\ge1$, consider the bounded homomorphism of $K$-modules
 $$
V[[z]]\ \buildrel {\psi_m}\over\to\ V^m\oplus V[[z]]: f(z)\mapsto (f(0),(\partial_z f)(0),\ldots,(\partial_z^{(m-1)})(0); \partial^{(m)}_z f(z))\ . 
 $$
 
\begin{lem}\label{lem-diffequ}
Assume that $V$ is torsion free as an abelian group. Then
\begin{itemize}
	\item [(i)]
the homomorphism $\psi_m$ is injective;
    \item [(ii)]
if in addition, $K$ is a $\bfQ$-algebra and the norm on $K$ induces the trivial norm on $\bfQ$, then $\psi_m$ is an isometric isomorphism.
\end{itemize}
\end{lem}

\begin{proof}
If $f(z) = \sum_{n=0}^\infty u_n z^n$, then
 $$
\psi_m(f(z)) = (u_0,u_1,\ldots,u_{m-1}; \sum_{n=m}^\infty \binom{n}{m} u_n z^{n-m})\ .
 $$
Thus, if $g(z)=\sum_{n=0}^\infty v_n z^n$, then $u_n=v_n$ for $0\le n\le m-1$ and $\binom{n}{m} (u_n-v_n)=0$ for $n\ge m$. Then the assumption in (i) implies that $u_n=v_n$ for $n\ge m$. The additional assumption in (ii) implies that $||\alpha v||=||v||$ for all $v\in V$ and all nonzero $\alpha\in\bfQ$, and this gives the required fact. 
\end{proof}

\subsection{Delta-function} \label{subsec-delta}

The formal {\it delta-function}
$\delta(z-w)$ is the element of $K[[z^{\pm1},w^{\pm1}]]$ defined
by $\delta(z-w) = \Delta({1\over{z-w}})$, i.e.,
$$
\delta(z-w) = \sum_{n\in\bfZ} z^n w^{-n-1}\ .
$$
Since $(z-w)\cdot{1\over{z-w}}=1$, it follows that
$(z-w)\delta(z-w)=0$. More generally, for each $i\ge0$ one has
$$
\partial_w^{(i)} \delta(z-w) =
\Delta\left({1\over{(z-w)^{i+1}}}\right) = \sum_{n\in\bfZ} \binom{n}{i} w^{n-i}z^{-n-1}
$$
and $(z-w)^{i+1}\partial_w^{(i)} \delta(z-w) = 0$. One also has
$$
\Res_z((z-w)^i \partial_w^{(i)} \delta(z-w)) = \Res_z \left(\Delta\left(
{1\over{(z-w)^{j-i+1}}}\right)\right) =\delta_{i,j}\ .
$$

Let $V$ be a Banach $K$-module.

\begin{lem}\label{lem-repres}
If an element $f\in
V[[z^{\pm1},w^{\pm1}]]$ can be represented in the form $f=
\sum_{i=0}^\infty g_i(w)\partial_w^{(i)} \delta(z-w)$, where $g_i(w)$ are
elements of $V[[w^{\pm1}]]$ such that $g_i(w)\to0$ as
$i\to\infty$, then $g_i(w)= \Res_z((z-w)^i f(z,w))$ and  
$||f||=\max\limits_i ||g_i||$.	
\end{lem}

\begin{proof}
The formula for $g_i(w)$ follows from the equality before the formulation; in particular, such a representation of $f(z,w)$ is unique. Since
$||\partial_w^{(i)} \delta(z-w)||\le1$, then $||f(z,w)||\le \max\limits_i
||g_i(w)||$, and since $g_i(w)= \Res_z ((z-w)^i
f(z,w))$, it follows that $||g_i(w)||\le ||f(z,w)||$.
\end{proof}

Let $V[[w^{\pm1}]]\{\delta(z-w)\}$ denote the closed $K$-submodule
of $V[[z^{\pm1},w^{\pm1}]]$ consisting of the elements that can be
represented in the form of Lemma \ref{lem-repres}. Then the lemma implies that such a representation is unique. Let also $V((w))\{\delta(z-w)\}$ denote the subspace of $V[[w^{\pm1}]]\{\delta(z-w)\}$ consisting of the elements for
which $g_i(w)\in V((w))$ for all $i\ge0$.

\begin{cor}\label{cor-repres}
For the operator $\Delta:
V((z,w))\{(z-w)^{-1}\} \to V[[z^{\pm1},w^{\pm1}]]$, one has
$\Ker(\Delta)=V((z,w))$ and ${\rm Im}(\Delta) =
V((w))\{\delta(z-w)\}$.
 \hfill$\square$	
\end{cor}

\begin{lem}\label{lem-repdelta}
(Decomposition theorem) The following properties of an element
$f(z,w)\in V[[z^{\pm1},w^{\pm1}]]$ are equivalent:
\begin{itemize}
	\item [(a)]
$f\in V[[w^{\pm1}]]\{\delta(z-w)\}$;
    \item[(b)]
$(z-w)^N f(z,w)\to0$ as $N\to\infty$.
\end{itemize}
\end{lem}

\begin{proof}
(a)$\Longrightarrow$(b) If $f$ is represented in the
form of (a), then
$$
(z-w)^N f(z,w) = \sum_{i=N}^\infty g_i(w) (z-w)^N\partial_w^{(i)} \delta(z-w)
$$
and, therefore, $||(z-w)^N f(z,w)|| \le \max\limits_{i\ge N}
||g_i|| \to 0$ as $N\to\infty$.

(b)$\Longrightarrow$(a) For $i\ge0$, we set $g_i(w)= \Res_z
((z-w)^i f(z,w))$. The assumption implies that $g_i\to0$ as
$i\to\infty$. We can therefore replace $f(z,w)$ by $f(z,w) -
\sum_{i=0}^\infty g_i(w)\partial_w^{(i)} \delta(z-w)$ and assume that $\Res_z
((z-w)^i f(z,w))=0$ for all $i\ge0$. In this case, {\it we claim
	that $f(z,w)=0$.} Indeed, the latter easily implies that
$f(z,w)\in V[[w^{\pm1}]][[z]]$ and, therefore, the claim follows
from Lemma \ref{lem-zminw}(i).
\end{proof}

\begin{rem}
	Since the operator of multiplication by $z-w$ is
	isometric on $V[[w^{\pm1}]]((z))$ (Lemma \ref{lem-zminw}(i)), a nonzero
	element of $V[[z^{\pm1},w^{\pm1}]]$ that possesses the equivalent
	properties of Lemma \ref{lem-repdelta}, cannot lie in $V[[w^{\pm1}]]((z))$.
\end{rem}

\subsection{Mutual locality of formal distributions} \label{subsec-localdist}

\begin{defin}\label{def-Lie}
A {\it Banach Lie $K$-algebra} is a Banach $K$-module $\gothg$ provided with the structure of a Lie $K$-algebra with the property $||[u,v]|| \le ||u||\cdot||v||$ for all $u,v\in \gothg$.
\end{defin}

For example, for any Banach $K$-module $V$, $\End(V)$ is a Banach Lie $K$-algebra. Let $\gothg$ be a Banach Lie $K$-algebra.
	
\begin{defin}\label{def-localdist}
Formal distributions $a(z),b(z)\in \gothg[[z^{\pm1}]]$ are said to be {\it mutually local} (or just {\it local}) if 	
 $$
(z-w)^N [a(z),b(w)] \to0\ {\it as\ } N\to\infty\ . 
 $$
\end{defin}

\begin{exam}\label{exam-localconst}
	Suppose that formal distributions $a(z)$ and $b(z)$ do not depend on $z$, i.e., they are just elements $a,b\in\gothg$. Then $a$ and $b$ are mutually local if and only if they commute, i.e., $[a,b]=0$.
\end{exam}

It is convenient to express the decomposition in $z$ of a formal distribution $a(z)\in \gothg[[z^{\pm1}]]$ in the form $a(z)=\sum_{n\in\bfZ} (a(z))_{(n)} z^{-n-1}$. For brevity, one also writes $a_{(n)} = (a(z))_{(n)} \in \gothg$, i.e., 
 $$
a(z)=\sum_{n\in\bfZ} a_{(n)} z^{-n-1}\ .
 $$

\begin{lem}\label{lem-localrep}
For every pair of mutually local formal distributions $a(z),b(z)\in \gothg[[z^{\pm1}]]$, the following is true:
\begin{itemize}
	\item[(i)]
there is a unique decomposition
$$
[a(z),b(w)]\ =\ \sum_{i=0}^\infty (a(w)_{(i)} b(w)) \partial^{(i)}_z \delta(z-w)\ ,
$$
where $a(w)_{(i)} b(w) \in\gothg[[z^{\pm1}]]$; moreover, one has
 $$
a(w)_{(i)} b(w)= \Res_z((z-w)^i [a(z),b(w)])\ ,
 $$ 
$a(w)_{(i)} b(w)\to0$ as $i\to\infty$ and
$||[a(z),b(w)]||=\max\limits_{i\ge0} \{||a(w)_{(i)} b(w)||\}$;
\item[(ii)]
for every $m,n\in\bfZ$, one has
 $$
[a_{(m)},b_{(n)}]\ =\ \sum_{ i=0}^\infty \binom{m}{i} (a(z)_{(i)} b(z))_{(m+n-i)} \ ;
 $$
 \item[(iii)]
for every $m\ge0$, the formal distributions $\partial_z^{(m)} a(z)$ and $b(z)$ are mutually local.
\end{itemize}
\end{lem}

\begin{proof}
The statement (i) follows from the definition and Lemmas \ref{lem-repres} and \ref{lem-repdelta}, and the statement (ii) follows from (i).

(iii) We prove the statement by induction on $m$. If $m=0$, it is true by the assumption. Assume that $m\ge1$ and the statement is true for all smaller values. Since multiplication by $z-w$ in $\gothg[[z^{\pm1},w^{\pm1}]]$ does not increase the norm, it suffices to show that, for every $\varepsilon>0$, there exists a sufficiently large $N$ such that the norm of the element
$(z-w)^N[\partial_z^{(m)} a(z),b(w)]$ is at most $\varepsilon$. Let $n$ be such that for all $0\le i\le m-1$ and $j\ge n$ the norm of the element $(z-w)^j
[\partial_z^{(i)} a(z),b(w)]$ is at most $\varepsilon$. One has
\begin{eqnarray}
	\partial_z^{(m)}((z-w)^N[a(z),b(w)]) &=& (z-w)^N
	[\partial^{(m)}_za(z),b(w)] {}\nonumber\\ & & +\sum_{i=1}^m \binom{N}{i} (z-w)^{N-i} [\partial_z^{(m-i)} a(z),b(w)]\ .
	{}\nonumber	
\end{eqnarray}
If $N\ge m+n$ then, by the induction, the left hand side and the second summand on the right hand side are at most $\varepsilon$. It follows that the same is true for the first summand.
\end{proof}

\subsection{Quantum fields and normally ordered products} \label{subsec-fields}
Let $V$ be a
Banach $K$-module. 

\begin{defin}\label{def-qfield}
A formal distribution
$$
a(z)=\sum_{n\in \bfZ} a_{(n)} z^{-n-1} \in \End(V)[[z^{\pm1}]]
$$
is called an {\it $\End(V)$-valued quantum field} (or just a {\it field}) if, for any $v\in V$, one has
$a(z)v\in V((z))$ (i.e., $a_{(n)}(v)\to0$ as $n\to+\infty$). 
\end{defin}

For example, all elements of $\End(V)((z))$ are fields and, in particular, the identity operator $I_V\in\End(V)$ is a field. It is easy to see that the space $\End(V) \langle\langle z
\rangle\rangle$ of all fields is closed in $\End(V)[[z^{\pm1}]]$ and is preserved by the operators $\partial^{(n)}_z$, $n\ge0$. 

\begin{lem}\label{lem-gnfield}
If $a(z)$ and $b(z)$ are fields, then for every $n\ge0$, the element $a(w)_{(n)} b(w) = {\rm Res}_z((z-w)^n [a(z),b(w)])$ (see Lemma \ref{lem-localrep}(i)) is a field.	
\end{lem}

\begin{proof}
If $a(w)_{(n)} b(w)=\sum_{p\in\bfZ}c_{(p)} w^{-p-1}$, one has
$$
c_{(p)} = \sum_{i=0}^n (-1)^{n-i} \binom{n}{i} [a_{(i)},b_{(p+n-i)}]\ .
$$
It follows that for every vector $v\in V$ one has
$$
||c_{(p)}(v)|| \le \max\limits_{0\le i\le n} \{||a_{(i)}||\cdot
||b_{(p+n-i)}(v)||, ||b_{(p+n-i)}(a_{(i)}(v))||\}\ .
$$
The latter number evidently tends to zero as $p$ tends to $+\infty$.
\end{proof}

For an element $a(z)\in
\End(V)[[z^{\pm1}]]$ as above, we set
$$
a(z)_+=\sum_{n\le -1} a_{(n)} z^{-n-1}\ {\rm and}\ a(z)_-=
\sum_{n\ge0} a_{(n)} z^{-n-1}\ .
$$
The {\it normally ordered product} of two fields $a(z),b(z)\in \End(V)\langle\langle z\rangle\rangle$ is defined by
$$
{\rm :}a(z)b(z){\rm :}\  = a(z)_+ b(z) + b(z) a(z)_-\ .
$$
For example, for any field $a(z)$, one has ${\rm :}a(z)I_V{\rm :}={\rm :}I_V a(z){\rm :}=a(z)$. 

\begin{lem}\label{lem-normprod}
The normally ordered product {\rm :}$a(z)b(z)${\rm :} is a well defined element of $\End(V)\langle\langle z\rangle\rangle$ whose norm is at most $||a(z)||\cdot||b(z)||$.
\end{lem}

\begin{proof}
One has
 $$a(z)_+ b(z) = \sum_{p\in\bfZ} c_p
z^{-p-1}\ {\rm and\ }b(z) a(z)_- = \sum_{p\in\bfZ} d_p z^{-p-1}\ ,
 $$ 
where
$$
c_p = \sum_{\buildrel {m+n=p-1}\over {m\le -1}} a_{(m)} b_{(n)}\ \ {\rm
	and}\ \ d_p = \sum_{\buildrel {m+n=p-1}\over {m\ge0}} b_{(n)} a_{(m)}\ .
$$
Although both infinite sums do not convergent with respect to the
Banach norm on $\End(V)$, they do convergent to bounded linear
operators in the weak topology. Let us show this, for example, for
the operator $c_p$. Given a vector $v\in V$, one has $||(a_{(m)} b_{(n)})(v)||\le ||a_{(m)}||\cdot ||b_{(n)}(v)||$. Since $m+n=p-1$ and $m\le
-1$, the latter number tends to zero as $n\to+\infty$, and so $c_p(v)$ is well
defined. Furthermore, given $\varepsilon>0$, there exists $n_0$
such that $||b_{(n)}(v)||\le \varepsilon$ for all $n\ge n_0$. It
follows that
$$
||c_p(v)|| \le \max\{\max\limits_{p\le n<n_0} ||a_{(p-1-n)}||\cdot
||b_{(n)}||\cdot||v||, \max\limits_{n\ge n_0}
||a_{(p-1-n)}||\varepsilon\}
$$
and, therefore, $||c_p||\le ||a(z)||\cdot||b(z)||$.
\end{proof}

\begin{lem}\label{lem-partial}
	For every $m\ge0$, one has
	$$
	\partial_z^{(m)} ({\rm :}a(z)b(z){\rm :}) = \sum_{i=0}^m  {\rm :}\partial_z^{(i)} a(z) \partial^{(m-i)}_z b(z){\rm :}\ .
	$$
\end{lem}

\begin{proof}
	The statement easily follows from the fact that the operator $\partial_z^{(m)}$ commutes with the operators $a(z)\mapsto a(z)_+$ and $a(z)\mapsto a(z)_-$.
\end{proof}

\begin{lem}\label{lem-tensorfield}
Let $U$ and $V$ be Banach $K$-modules, and let $a(z)\in \End(U)\langle\langle z\rangle\rangle$ and $b(z)\in \End(V)\langle\langle z\rangle\rangle$. Then
\begin{itemize}
    \item[(i)]
for every $n\in\bfZ$, the sum $\sum_{i+j=n-1} a_{(i)} \otimes b_{(j)}$ converges in the week topology of $U\hatotimes_K V$ to a bounded operator $c_{(n)}\in\End(U\hatotimes_K V)$ of norm at most $||a(z)||\cdot ||b(z)||$;
    \item[(ii)]
for every $x\in U\hatotimes_K V$, $c_{(n)}(x) \to0$ as $n\to+\infty$ and, therefore, the sum $\sum_{n\in\bfZ} c_{(n)} z^{-n-1}$ is an element of $\End(U\hatotimes_K V)\langle\langle z\rangle\rangle$, denoted by $a(z)\hatotimes b(z)$;
    \item[(iii)]
if $a'(z),a''(z)\in\End(U)\langle\langle z\rangle\rangle$ and $b'(z),b''(z)\in\End(V)\langle\langle z\rangle\rangle$ are pairs of mutually local fields, then the fields $a'(z)\hatotimes b'(z)$ and $a''(z) \hatotimes b''(z)$ in $\End(U\hatotimes V)\langle\langle z\rangle\rangle$ are mutually local.
\end{itemize}
\end{lem}

\begin{proof}
(i) Let $x\otimes y$ be a nonzero element of $U\otimes_K V$. Since $a(z)$ and $b(z)$ are fields, it follows that, for any $\varepsilon>0$, there exists $N\ge1$ such that, for all $i,j\ge N$, one has $||a_{(i)}(x)||\le \frac{\varepsilon}{||b(z)||\cdot||y||}$ and $||b_{(j)}(y)||\le \frac{\varepsilon}{||a(z)||\cdot||x||}$. Thus, if $i\ge N$, one has 
 $$
||(a_{(i)}\otimes b_{(j)})(x\otimes y)||\le ||a_{(i)}(x)||\cdot||b_{(j)}(y)|| \le \frac{\varepsilon}{||b(z)||\cdot||y||} \cdot  ||b_{(j)}(y)|| \le \varepsilon\ ,
$$
and the similar inequality holds for $j\ge N$. This implies that the sum considered is weakly converges and induces a bounded operator $c_{(n)}$ on $U\hatotimes_K V$ of norm at most $||a(z)||\cdot ||b(z)||$, i.e., (i) is true. 

(ii) If $n\ge 2N+1$ and $i+j=n-1$, then either $i$ or $j$ is at least $N$, and so $||c_{(n)}(x\otimes y)||\le\varepsilon$. This implies (ii).

(iii) The commutator $[a'(z)\hatotimes b'(z),a''(w) \hatotimes b''(w)]$, considered as an element of $\End(U\hatotimes V) [[z^{\pm 1},w^{\pm 1}]]$, is the image of the following element 
$$
[a'(z),a''(w)] \otimes b'(z)b''(w) + a''(w)a'(z) \otimes [b'(z),b''(w)]\ .
$$
Since the norms of the elements $b'(z)b''(w)$ and $a''(w)a'(z)$ are bounded by a constant, we see that the norm of the commutator, multiplied by $(z-w)^N$, tends to zero as $N$ tends to infinity.
\end{proof}

Lemma \ref{lem-tensorfield} implies that, for Banach $K$-modules $U$ and $V$, there is a well defined bounded homomorphism of Banach $K$-modules
$$
 \End(U)\langle\langle z\rangle\rangle \hatotimes \End(V)\langle\langle z\rangle\rangle \to \End(U\hatotimes V)  \langle\langle z\rangle\rangle: a(z)\otimes b(z) \mapsto a(z)\hatotimes b(z)\ ,
 $$
and one has $||a(z)\hatotimes b(z)||\le ||a(z)||\cdot ||b(z)||$.

\subsection{$n$-th products of quantum fields}
Let $V$ be a Banach $K$-module. By Lemma \ref{lem-gnfield}, for each $n\ge0$, the $n$-th product $a(z)_{(n)} b(z)$ of fields is a field. One can extend these $K$-bilinear operations on the space of fields to arbitrary $n\in\bfZ$. Namely, if $n=-m-1<0$, we set
$$
a(z)_{(n)} b(z) =\  {\rm :}(\partial_z^{(m)} a(z)) b(z){\rm :}\ .
$$
(The latter is a field, by Lemma \ref{lem-normprod} and the fact that the space
$\End(V)\langle\langle z\rangle\rangle$ is preserved by the
operators $\partial_z^{(m)}$, $m\ge0$.) For example, for any field $a(z)$, one has $a(z)_{(n)} I_V=0$ for $n\ge0$, and $a(z)_{(n)} I_V=\partial^{(m)}_z a(z)$, for $n=-m-1<0$. Notice that one always has $||a(z)_{(n)} b(z)||\le ||a(z)||\cdot||b(z)||$.

\begin{exam}\label{exam-nprodconst}
As in Example \ref{exam-localconst}, suppose that fields $a(z)$ and $b(z)$ do not depend on $z$, i.e., they are just operators $a,b\in\End(V)$. Then their $n$-th product is the composition $ab$, if $n=-1$, and zero, if $n\not=0$.
\end{exam}

\begin{lem}\label{lem-prodform}
For every $n\in\bfZ$, one has
 $$
a(w)_{(n)} b(w) =\Res_z(a(z)b(w)i_{z,w}((z-w)^n) - b(w)a(z)i_{w,z}((z-w)^n))\ . 
 $$
\end{lem}

\begin{proof}
By the definition, if $n\ge0$, then $a(w)_{(n)} b(w)=g_n(w)$ for $g_n(w) = \Res_z((z-w)^n [a(z),b(w)])$. This coincides with the right hand side of the required formula. If $n=-m-1<0$, the formula follows from the following equalities:
 $$
\Res_z\left(a(z) i_{z,w}\left({1\over {(z-w)^{m+1}}}\right)\right) =  \partial^{(m)}_w a(w)_+\ {\rm and}
 $$
 $$
\Res_z\left(a(z) i_{w,z}\left({1\over {(z-w)^{m+1}}}\right)\right) = \partial^{(m)}_w a(w)_-\ , 
 $$
which easily follow from the formulas for $i_{z,w}\left({1\over {(z-w)^{m+1}}}\right)$ and $i_{w,z}\left({1\over {(z-w)^{m+1}}}\right)$ established in \S \ref{subsec-deres}. 
\end{proof}

\begin{cor}\label{cor-preBorcherds}
	If fields $a(z)$ and $b(z)$ are mutually local, then for each $n\in\bfZ$, one has
	$$
	a(z)b(w)i_{z,w}((z-w)^n) - b(w)a(z)i_{w,z}((z-w)^n) = \sum_{j=0}^\infty (a(w)_{(n+j)} b(w)) \partial^{(j)}_w \delta(z-w)\ . 
	$$
\end{cor}

\begin{proof}
If an integer $N$ is large enough so that $N+n\ge0$, then the element on the left hand side multiplied by $(z-w)^N$ is equal to $(z-w)^{N+n}[a(z),b(w)]$. Since $a(z)$ and $b(z)$ are mutually local, the latter tends to zero as $N\to+\infty$. This means that the element on the left hand side satisfy the condition (b) of Lemma \ref{lem-repdelta} and, therefore, this element is of the form $\sum_{j=0}^\infty g_j(w) \partial^{(j)}_w \delta(z-w)$ with $g_j(w)\to0$ as $j\to\infty$. By Lemma \ref{lem-repres}, the element $g_j(w)$ is equal to
\begin{eqnarray}
& & \Res_z((a(z)b(w)i_{z,w}((z-w)^n) - b(w)a(z)i_{w,z}((z-w)^n))(z-w)^j) 	{}\nonumber	\\ 
&=& \Res_z(a(z)b(w)i_{z,w}((z-w)^{n+j}) - b(w)a(z)i_{w,z}((z-w)^{n+j})) {}\nonumber	\\
&=& a(w)_{(n+j)} b(w)\ . 
	{}\nonumber	
\end{eqnarray}
The required statement follows.
\end{proof}

\begin{lem}\label{lem-partialn}
For any pair of fields $a(z)$ and $b(z)$ and every $n\in\bfZ$, the following is true:
\begin{itemize}
	\item [(i)]
 $
\partial_z(a(z)_{(n)} b(z)) = (\partial_z a(z))_{(n)} b(z) + a(z)_{(n)} \partial_z b(z)\ ;
 $
    \item[(ii)]
if $T$ is a bounded $K$-linear endomorphism of $V$ such that $[T,a(z)] = \partial_z a(z)$ and $[T,b(z)] = \partial_z b(z)$, then
 $$
[T,a(z)_{(n)} b(z)] = T a(z)_{(n)} b(z) + a(z)_{(n)} T b(z)\ . 
 $$
\end{itemize}
\end{lem}

\begin{proof}
(i) If $n<0$, the equality follows from the definition and the similar equality for normally ordered products from Lemma \ref{lem-partial} (for $m=1$). Suppose $n\ge0$. Then $a(w)_{(n)}b(w) = \Res_z((z-w)^n [a(z),b(w)])$, and $\partial^{(m)}_w(a(w)_{(n)}b(w))$ is equal to
\begin{eqnarray}
	& & \Res_z(a(z) \partial^{(m)}_w((z-w)^n b(w)) - \partial^{(m)}_w((z-w)^n b(w) a(z))) {}\nonumber\\ &=& \Res_z(a(z) \partial_w((z-w)^n) b(w) - \partial_w((z-w)^n) b(w) a(z))  {}\nonumber\\ &+& \Res_z(a(z) (z-w)^n \partial_w b(w) - (z-w)^n \partial_w b(w) a(z))\ .
	{}\nonumber
\end{eqnarray}
Substituting the equality $\partial_w((z-w)^n) = -\partial_z((z-w)^n)$ and then the equality $\Res_z(a(z) \partial_z((z-w)^n)) = - \Res_z(\partial_z a(z) (z-w)^n)$, we get the required fact.
 \medskip
 
(ii) By Lemma \ref{lem-prodform}, the left hand side is equal to
\begin{eqnarray}
	& & \Res_z(Ta(z)b(w)i_{z,w}((z-w)^n) - a(z)b(w)T i_{z,w}((z-w)^n) {}\nonumber\\
	&-& Tb(w)a(z) i_{w,z}((z-w)^n) - b(w)a(z)T i_{w,z}((z-w)^n)) {}\nonumber\\
	&=& \Res_z (Ta(z)b(w)i_{z,w}((z-w)^n) - a(z)Tb(w)i_{z,w}((z-w)^n) {}\nonumber\\
	&+& a(z)Tb(w)i_{z,w}((z-w)^n) - a(z)b(w)T i_{z,w}((z-w)^n)  {}\nonumber\\	
	&-& Tb(w)a(z) i_{w,z}((z-w)^n) + b(w)T a(z) i_{w,z}((z-w)^n) {}\nonumber\\	
	&-& b(w)Ta(z) i_{w,z}((z-w)^n) - b(w)a(z)T i_{w,z}((z-w)^n)) {}\nonumber\\	
	&=& \Res_z ([T,a(z)] b(w) i_{z,w}((z-w)^n) - b(w) [T,a(z)] i_{w,z}((z-w)^n)) {}\nonumber\\	
	&+& \Res_z (a(z)[T,b(w)] i_{z,w}((z-w)^n) - [T,b(w)] a(z) i_{w,z}((z-w)^n)).
	{}\nonumber
\end{eqnarray}
The assumption implies that the latter expression is equal to the right hand side of the required equality.
\end{proof}

A field $a(z)$ is said to be {\it translation covariant with respect} to an endomorphism $T\in\End(V)$ if $[T,a(z)] = \partial_z a(z)$.

\begin{cor}\label{cor-translation}
If fields $a(z)$ and $b(z)$ are translation covariant with respect to an endomorphism $T\in\End(V)$, then so are their $n$-th products $a(z)_{(n)} b(z)$ for all $n\in\bfZ$.
 \hfill$\square$
\end{cor}

\begin{lem}\label{lem-vacuum}
Assume that for $a(z)$ and $b(z)$ and for some element $v\in V$, one has $a(z)v,b(z)v\in V[[z]]$. If $c(z)=a(z)_{(n)}b(z)$ for $n\in\bfZ$, then
\begin{itemize}
\item[(i)]
$c(z)v \in V[[z]]$;
\item[(ii)]
$c_{(-1)}v = a_{(n)} (b_{(-1)}v)$.
\end{itemize}
\end{lem}

\begin{proof}
(i) By Lemma \ref{lem-prodform},  $c(w)v$ is equal to
 $$
c(w)v = \Res_z(a(z)b(w)i_{z,w}((z-w)^n))v - \Res_z(b(w)a(z)i_{w,z}((z-w)^n))v \ .
 $$
Since $a(z)v\in V[[z]]$ and the series $i_{w,z}((z-w)^n$ has no negative powers of $z$, the second summand is equal to zero. Since $b(w)v\in V[[w]]$ and the element $i_{z,w}((z-w)^n$ has no negative powers of $w$, the first summand lies in $V[[w]]$. The statement follows.

(ii) We can substitute the value $w=0$ to the above equality, and we get $c(0)v = \Res_z(a(z)b(0)z^n)v=a_{(n)} (b(0)v)$. Since $c(z),b(z)\in V[[z]]$, one has $c(0)=c_{(-1)}$ and $b(0)=b_{(-1)}$, and the required equality follows. 
\end{proof}

\begin{lem}\label{lem-Dong}
(Dong's Lemma) 
If $a(z)$, $b(z)$ and $c(z)$ are pairwise mutually local fields, then for each $n\in\bfZ$ the fields $a(z)_{(n)} b(z)$ and $c(z)$ are mutually local.
\end{lem}

\begin{proof}
By Lemma \ref{lem-prodform}, it suffices to show that the
element $(w-u)^N P$ tends to zero as $N\to\infty$, where
\begin{eqnarray}
	P &=& (a(z)b(w)i_{z,w}((z-w)^n) -b(w)a(z) i_{w,z} ((z-w)^n)) c(u) {}\nonumber\\ & & -\ c(u)(a(z)b(w)i_{z,w}((z-w)^n) -b(w)a(z) i_{w,z} ((z-w)^n)) {}\nonumber\\ &=& [a(z)b(w),c(u)] i_{z,w}((z-w)^n) - [b(w)a(z),c(u)] i_{w,z} ((z-w)^n)
{}\nonumber	
\end{eqnarray}
Given $\varepsilon>0$, let $k$ be such that $k\ge\max\{n,-n\}$ and for all $m\ge k$ the elements $(z-w)^m[a(z),b(w)]$, $(z-u)^m[a(z),c(u)]$ and
$(w-u)^m[b(w),c(u)]$ are of norm at most $\varepsilon$. For $N=3k$
we have
$$
(w-u)^{3k} = (w-u)^k \sum_{m=0}^{2k} \binom{2n}{m}
(w-z)^m (z-u)^{2k-m}\ .
$$
To verify the required fact, we evaluate separately summands in
$(w-u)^{3k}P$ for $0\le m\le k$ and for $k+1\le m\le 2k$.

(1) If $0\le m\le k$, we substitute to the formula for $P$ the equalities
 $$
[a(z)b(w),c(u)] = a(z)[b(w),c(u)] + [a(z),c(u)]b(w)\ {\rm and} 
 $$
 $$
[b(w)a(z),c(u)] = b(w)[a(z),c(u)] + [b(w),c(u)] a(z)\ , 
 $$
and deduce from the assumptions that the element
$(w-u)^k(z-u)^{2k-m} P$ is of norm at most
$\varepsilon\max\{||a(z)||,||b(w)||\}$.

(2) If $k+1\le m\le 2k$, we use the equality
$$
P= [[a(z),b(w)],c(u)] i_{z,w}((z-w)^n) + [b(w)a(z),c(u)]\Delta((z-w)^n)\ .
$$
The product of the first summand by $(z-w)^m$ is of norm at most $\varepsilon||c(u)||$. If $n\ge0$, then $\Delta((z-w)^n)=0$, and so the second summand is equal to zero. If $n=-p-1<0$, then $\Delta((z-w)^p)= \partial_w^{(p)} \delta(z-w)$. Since $(z-w)^{p+1} \partial_w^{(p)} \delta(z-w)=0$ (see \S \ref{subsec-delta}) and $k\ge-n=p+1$, it follows the second summand vanishes by multiplication by $(w-z)^m$. 
\end{proof}

\subsection{The closed $K$-submodule generated by mutually local fields}
Let $\calF$ be a family of mutually local fields. The {\it $K$-submodule generated by $\calF$}  is the minimal $K$-submodule $\calF_{\rm min}$ of $\End(V)\langle\langle z\rangle\rangle$, which contains all fields from $\calF\cup \{I_V\}$ and is preserved by all of the $n$-products. The {\it closed $K$-submodule generated by $\calF$} is the closure $\overline\calF_{\rm min}$ of $\calF_{\rm min}$ in $\End(V)\langle\langle z\rangle\rangle$.

\begin{lem}\label{lem-localclosed}
In the above situation, the following is true:
\begin{itemize}
	\item [(i)]
all fields from $\overline\calF_{\rm min}$ are mutually local;
    \item [(ii)]
$\overline\calF_{\rm min}$ is preserved by all of the $n$-products;
    \item [(iii)]
given $v\in V$, if $\varphi(z)v\in V[[z]]$ for all $\varphi(z)\in\calF$, then the same holds for all fields $\varphi(z)\in \overline\calF_{\rm min}$;
    \item [(iv)]
if all fields from $\calF$ are translation covariant with respect to an endomorphism $T\in\End(V)$, then the same holds for all fields $\varphi(z)\in \overline\calF_{\rm min}$.
\end{itemize}
\end{lem}

\begin{proof}
Dong's Lemma \ref{lem-Dong} implies mutual locality of all fields from $\calF_{\rm min}$. It is also clear that the property (ii) extend to $\calF_{\rm min}$ and the properties (iii)-(iv) extend to $\overline\calF_{\rm min}$. Thus, in order to prove (i) and (ii), we may replace $\calF$ by $\calF_{\rm min}$. 

Let $\{a^i(z)\}_{i\ge1}$ and $\{b^i(z)\}_{i\ge1}$ be sequences of fields from $\calF$ that  converge to fields $a(z)$ and $b(z)$, respectively. For (i) (resp. (ii)), we have to show that $a(z)$ and $b(z)$ are mutually local (resp. the sequence $a^i(z)_{(n)} b^i(z)$ converges to $a(z)_{(n)} b(z)$). For each $i\ge1$, one has
$$
[a(z),b(w)] - [a^i(z),b^i(w)] = [a(z)-a^i(z),b(w)]+ [a^i(z),b(w)-b^i(w)]
$$
$$
({\rm resp.\ } a(z)_{(n)}b(z) - a^i(z)_{(n)} b^i(z) = (a(z)-a^i(z))_{(n)} b(z) + a^i(z)_{(n)} (b(z)-b^i(z)))\ .
$$
Given $\varepsilon>0$, we can find $i\ge1$ such that $||a(z)-a^i(z)|| < \varepsilon$, $||b(z)-b^i(z)|| < \varepsilon$, and $||a^i(z)||=||a(z)||$. Then the norms of the two summands on the right hand side are less than $\varepsilon\max\{||a(z)||,||b(z)||\}$. This immediately implies (ii). Furthermore, we can find $N\ge1$ such that, for all $j\ge N$, one has $||(z-w)^N [a^i(z),b^i(w)]||< \varepsilon$. Multiplication by $(z-w)^j$ does not increase the norms of the above two summands, but the norm of the second summand on the left hand side, multiplied by $(z-w)^j$ is less than $\varepsilon$. It follows that $(z-w)^N [a(z),b(w)]||\to0$ as $N\to\infty$.
\end{proof}

\section{Vertex $K$-algebras}\label{sec-vertex}

\subsection{Definition of a non-Archimedean vertex algebra}\label{subsec-defvertex}
In this section, we assume that the commutative Banach ring $K$ as well as all Banach $K$-modules considered are torsion free as abelian groups.

\begin{defin}\label{def-vertex}
A {\it vertex $K$-algebra} is a collection of data:
\begin{itemize}
\item[(a)] 
({\it space of states}) a torsion free Banach $K$-module $V$;
\item[(b)]
({\it vacuum vector}) an element $|0\rangle \in V$;
\item[(c)]
({\it translation operator}) a bounded $K$-linear operator $T:V\to V$ with $T|0\rangle=0$;
\item[(d)]
({\it space of fields}) a closed $K$-submodule ${\rm Fld}(V) \subset \End(V)\langle\langle z\rangle\rangle$ that contains the identity operator $I_V\in\End(V)$;
\end{itemize}
with the following properties:
\begin{itemize}
\item[(V.1)]
$\varphi(z)|0\rangle \in V[[z]]$ for all $\varphi(z)\in {\rm Fld}(V)$;
\item[(V.2)]
all fields from ${\rm Fld}(V)$ are translation covariant with respect to $T$;
\item[(V.3)]
all fields from ${\rm Fld}(V)$ are mutually local;
\item[(V.4)]
${\rm Fld}(V)$ is preserved under all of the $n$-th products;
\item[(V.5)]
for the image $V'$ of the following bounded homomorphism
$$
{\rm Fld}(V)\to V: \varphi(z)\mapsto\varphi_{(-1)}|0\rangle\ ,
$$
its $K$-saturation $\{v\in V \big\vert \lambda v\in V'$ for a nonzero $\lambda\in K\}$ is dense in $V$.
\end{itemize}
 \end{defin}
 
For brevity, the above object is referred to as a vertex $K$-algebra $V$. The bounded homomorphism from (V.5) is denoted by $fs$ (field-statee). 
 
\begin{lem}\label{lem-fs}
The field-state homomorphism $fs:{\rm Fld}(V)\to V$ is injective.
\end{lem}

\begin{proof}
Let $\varphi(z)\in {\rm Fld}(V)$ be a field from the kernel of the homomorphism $fs$, and let $\psi(z)\in{\rm Fld}(V)$. Since $\varphi(z)$ and $\psi(z)$ are mutually local, it follows that
$$
(z-w)^N \varphi(z)\psi(w) |0\rangle - (z-w)^N \psi(w)\varphi(z) |0\rangle \to 0\ {\rm as\ } N\to+\infty\ .
$$
The property (V.2) means that $[T,\varphi_{(n)}]=(-n)\varphi_{(n-1)}$ and, since $T|0\rangle=0$, it follows that $T\varphi_{(n)}|0\rangle = (-n)\varphi_{(n-1)}|0\rangle$ for all $n\in\bfZ$. Furthermore, since $fs(\varphi(z))=\varphi_{(-1)}|0\rangle=0$ and $V$ is torsion free as an abelian group, it follows, by induction, that $\varphi_{(n)}|0\rangle=0$ for all $n<0$. Finally, by the property (V.1), $\varphi_{(n)}|0\rangle=0$ for all $n\ge0$. Thus, $\varphi(z)|0\rangle=0$ and, therefore, the second summand of the above expression is zero. It follows that $(z-w)^N \varphi(z)\psi(w) |0\rangle \to0$ as $N\to +\infty$. Finally, since $\psi(w)|0\rangle \in V[[w]]$, we can substitute $w=0$ in the left hand side, and we get $z^N \varphi(z) \psi_{(-1)}|0\rangle \to0$ as $N\to +\infty$. Since the norm of the element considered is equal to the norm of $\varphi(z) \psi_{(-1)}|0\rangle$, it follows that the latter element is equal to zero. Since $V$ is torsion free as a $K$-module, the property (V.5) then implies that $\varphi(z)=0$.
\end{proof}

Lemma \ref{lem-fs} implies that there is a bijective homomorphism of $K$-modules
 $$
V'\to{\rm Fld}(V): a\mapsto Y(a,z) =   \sum_{n\in\bfZ} a_{(n)} z^{-n-1}\ ,
 $$
which is inverse to $fs$ and denoted by $sf$ (state-field). It can be viewed as an unbounded operator $sf:V\to {\rm Fld}(V)$ with the domain $V'$. The $K$-module $V'$ is complete with respect to the Banach norm, defined by $||a||' = ||Y(a,z)||$, the homomorphism $sf:(V',||\ ||') \to {\rm Fld}(V)$ is isometric, and the canonical embedding $(V',||\ ||')\to V$ is bounded.

\begin{defin}\label{def-admiss}
A vertex $K$-algebra $V$ is said to be {\it admissible} if the field-state homomorphism of Banach $K$-modules $fs:{\rm Fld}(V)\to V$ is admissible.
\end{defin}

In this case, $V'$ is a closed $K$-submodule of $V$ and the state-field operator $sf: V'\to{\rm Fld}(V)$ is bounded. If $K$ is a non-Archimedean field with nontrivial valuation, then by the Banach open mapping theorem, a vertex $K$-algebra $V$ is admissible if and only if $V'=V$.

\begin{exam}\label{exam-T0}
Let $V$ be a vertex $K$-algebra, whose translation operator $T$ is zero. Then all fields in ${\rm Fld}(V)$ do not depend on $z$ and, by Examples \ref{exam-localconst} and \ref{exam-nprodconst}, ${\rm Fld}(V)$ is a commutative Banach $K$-subalgebra of $\End(V)$ and, by the property (V.5), the $K$-saturation of the image $V'$ of the field-state homomorphism ${\rm Fld}(V) \to V: \varphi\mapsto fs(\varphi)=\varphi|0\rangle$ is dense in $V$. Conversely, given a torsion free Banach $K$-module $V$ and an element $|0\rangle\in V$, any commutative Banach $K$-subalgebra of $\End(V)$ with the latter property defines the structure of a vertex $K$-algebra on $V$ with zero translation operator. For example, let $K$ be the field of $p$-adic numbers $\bfQ_p$, $V=\bfQ_p\{x\}$, and $|0\rangle = 1+px+p^2x^2+\ldots$. Let also ${\rm Fld}(V)$ be the commutative Banach $K$-subalgebra of $\End(V)$ consisting of the operators that correspond to infinite sequences $\varphi=(\lambda_0,\lambda_1,\ldots)$ of elements of $\bfQ_p$ with $||\varphi||=\sup_n\{|\lambda_n|\}<\infty$ and act on $V$ as follows: $\varphi (\sum_{n=0}^\infty \alpha_n x^n) = \sum_{n=0}^\infty \lambda_n\alpha_n x^n$. Then $V$ is a vertex $K$-algebra, which is not admissible.  Indeed, for $\varphi$ as above, one has $fs(\varphi) = \sum_{n=0}^\infty \lambda_n p^n x^n$. Thus, if for $n\ge0$, $\varphi_n$ is the operator defined by $\varphi_n(x^m) = \delta_{m,n} x^m$, one has $||\varphi_n||=1$ and $||fs(\varphi_n)||=|p|^n$. The latter tends to zero when $n$ tends to infinity, i.e., the operator $fs$ is not admissible.
\end{exam}

Let again $V$ be an arbitrary vertex $K$-algebra. It follows from the definition of the operator $sf$ that, for any $a\in V'$, one has $a=a_{(-1)}|0\rangle$. In particular, if $|0\rangle=0$, then $V=0$. If $|0\rangle\not=0$, then $||a||\le |||0\rangle||\cdot ||a||'$ and $|||0\rangle||'=1$. Furthermore, by the property (V.2), $T a_{(n)}|0\rangle = -na_{(n-1)}|0\rangle$ and, therefore, $Ta=a_{(-2)}|0\rangle$.

\begin{lem}\label{lem-Tn}
For every $n\ge0$, the following is true:
\begin{itemize}
\item[(i)]
there is a well defined $K$-linear operator $\frac{1}{n!} T^n: V'\to V'$, and one has $\frac{1}{n!} T^n=T^{(n)}$ for the $K$-linear operator $T^{(n)}: V'\to V'$, defined by $T^{(n)} a= a_{(-n-1)}|0\rangle$;
\item[(ii)]	
there is a well defined $K$-linear operator $e^{zT}: V'\to V'[[z]]$, defined by $e^{zT} a = \sum_{n=0}^\infty \frac{1}{n!} T^n a$, and one has $e^{zT}a = Y(a,z) |0\rangle$;
\item[(iii)]
if the operators of multiplication by integers on $V$ are admissible , then the operators $\frac{1}{n!} T^n$ are extended to bounded operators $\oV'\to \oV{}'$ on the closure $\oV{}'$ of $V'$ in $V$.
\end{itemize}
\end{lem}

\begin{proof}
The equality $T a_{(n)}|0\rangle = -na_{(n-1)}|0\rangle$ implies that
 $$
T^m a_{(n)}|0\rangle = (-1)^m n(n-1)\cdot\ldots\cdot (n-m+1) a_{(n-m)}|0\rangle
 $$
for all $m\ge0$ and, in particular, $T^m a_{(-1)}|0\rangle = m! a_{(-m-1)} |0\rangle$. Since $a_{(-1)}|0\rangle=a$ and $V$ is torsion free as an abelian group, it follows that $\frac{1}{m!} T^m a = a_{(-m-1)}|0\rangle$. This implies (i). Since $||T^{(m)} a|| \le ||a||'\cdot |||0\rangle||$, the element $e^{zT}a$ lies in $V'[[z]]$, and the statement (ii) follows. Under the assumption of (iii), there exists a constant $C>0$ such that $||v||\le C||m! v||$ for all $v\in V$. This implies that $||T^{(m)}||\le C ||T^m||$, and (iii) follows.
\end{proof}

\begin{rems}\label{rem-equivnorm}
(i) The structure of a vertex $K$-algebra on $V$ defines, for every $n\in\bfZ$, a $K$-bilinear operation $V'\times V\to V: (a,b)\mapsto a_{(n)} b$. It follows from the definition that, for each pair $(a,b)\in V'\times V$, one has $||a_{(n)}b|| \le ||a||'\cdot||b||$ and $a_{(n)} b\to0$ as $n\to+\infty$. Recall that, if $V$ is admissible, then $V'$ is closed in $V$ and the norms $||\ ||'$ and $||\ ||$ on it are equivalent.

(ii) If $K$ is a field, the $K$-saturation of $V'$ coincides with $V'$.

(iii) If the Banach norms on $K$ and $V$ are trivial, then any vertex $K$-algebra is admissible but $V'$ does not necessarily coincides with $V$ (see \S \ref{subsec-bosfer}).
\end{rems}

\begin{defin}\label{def-morphism}
	A {\it homomorphism of vertex $K$-algebras $U\to V$} is a homomorphism of Banach $K$-modules $\varphi:U\to V$ such that 
	\begin{itemize}
		\item [(H.1)]
		$\varphi(|0\rangle) = |0\rangle$;
		\item[(H.2)]
		$\varphi\circ T = T\circ \varphi$; 	
		\item[(H.3)]	
		$\varphi(U')\subset V'$;   
		\item[(H.4)]    
		$\varphi(a_{(n)} b) = \varphi(a)_{(n)} \varphi(b)$ for all $a\in U'$, $b\in U$, and $n\in \bfZ$;
		\item[(H.5)]
		the induced homomorphism ${\rm Fld}(U)\to{\rm Fld}(V): Y(a,z) \mapsto Y(\varphi(a),z)$ is bounded.
	\end{itemize}
\end{defin}

\begin{defin}\label{def-ideal}
	An {\it ideal of a vertex $K$-algebra $V$} is a closed $K$-submodule $J\subset V$ such that $V/J$ is a torsion free $K$-module,  $TJ\subset J$ and $\psi_{(n)}J\subset J$ for all $\psi(z)=\sum_{n\in\bfZ} \psi_{(n)} z^{-n-1}\in {\rm Fld}(V)$ and $n\in\bfZ$.
\end{defin}

\begin{lem}\label{lem-ideal}
	Given an ideal $J$ of a vertex $K$-algebra $V$, the quotient $V/J$ has a canonical structure of a vertex $K$-algebra for which the canonical map $\pi:V\to V/J$ is a homomorphism of vertex $K$-algebras.
\end{lem}

\begin{proof}
	The vertex $K$-algebra structure on $V/J$ is defined by the vacuum vector $\pi(|0\rangle)$, the translation operator $\pi(T)$, and the space of fields ${\rm Fld}(V/J)$, which is the closure of the $K$-module of $\End(V/J)$-valued fields of the form $\pi(\psi(z)) = \sum_{n\in\bfZ} \pi(\psi_{(n)})z^{-n-1}$ for $\psi(z) \in {\rm Fld}(V)$. (Here $\pi(T)$ and $\pi(\psi_{(n)})$ are the operators on $V/J$, induced by $T$ and $\psi_{(n)}$, respectively.) Indeed, the latter $K$-module obviously possesses the properties (V.1)-(V.5) and, by Lemma \ref{lem-localclosed}, the same properties hold for its closure ${\rm Fld}(V/J)$. This implies the claim.
\end{proof} 

Lemma \ref{lem-ideal} and injectivity of the state-field homomorphisms of vertex $K$-algebras imply that, for each element $a\in J\cap V'$, one has $a_{(n)} V\subset J$. Notice also that an ideal $J$ of $V$ is nontrivial (i.e., $J\not=V$) if and only if $|0\rangle\not\in J$.

\subsection{Basic properties of vertex algebras}
Let $V$ be a vertex $K$-algebra. 

\begin{lem}\label{lem-Goddard}
	(Goddard's Uniqueness Theorem) Suppose that a field $X(z)$ is mutually local with respect to all fields from ${\rm Fld}(V)$ and such that $X(z)|0\rangle = \varphi(z)|0\rangle$ for some $\varphi(z)\in{\rm Fld}(V)$. Then $X(z)= \varphi(z)$.
\end{lem}

\begin{proof}
	First of all, replacing $X(z)$ by $X(z)-\varphi(z)$, we may assume that $X(z)|0\rangle=0$, and our purpose is to show that
	in this case $X(z)=0$. By the assumption, for every $\psi(z)\in{\rm Fld}(V)$ and
	$N\ge0$, we have
	$$
	(z-w)^N[X(z),\psi(w)] |0\rangle = (z-w)^N X(z) \psi(w)  |0\rangle
	$$
	and, by the vacuum axiom, the right hand side is an element of
	$V[[z^{\pm 1},w]]$. By Lemma \ref{lem-zminw}(i), multiplication by $(z-w)^N$ on $V[[z^{\pm 1},w]]$ is an isometry and, therefore, one has
	$$
	||(z-w)^N[X(z),\psi(z)] |0\rangle|| = ||X(z) \psi(w)  |0\rangle|| \ge ||X(z)\psi_{(-1)}|0\rangle||\ .
	$$
	By the locality assumption, the left hand side tends to zero as $N\to\infty$. It follows that $X(z)b=0$ for all $b\in V'$. Since the $K$-saturation of $V'$ is dense in $V$, we get $X(z)=0$.
\end{proof}

\begin{lem}\label{lem-basic}
(i) The $K$-module ${\rm Fld}(V)$ is preserved by the operators $\partial_z^{(n)}$, and $\partial_z^{(n)} Y(a,z) = Y(T^{(n)} a,z)$ for all $a\in V'$ and $n\ge0$;

(ii) $Y(a,z)_{(n)}Y(b,z) = Y(a_{(n)}b,z)$ for all $a,b\in V'$ and $n\in\bfZ$.
\end{lem}

\begin{proof}
(i) One has $\partial^{(m)}_z Y(a,z) = Y(a,z)_{(n)} I_V$ for $n=-m-1<0$, and since ${\rm Fld}(V)$ contains $I_V$ and is preserved by the $n$-product, it follows that $\partial^{(m)}_z Y(a,z) \in {\rm Fld}(V)$. Furthermore, we apply Lemma \ref{lem-Goddard} to the field $X(z) = \partial_z^{(m)} Y(a,z)$. By Lemma \ref{lem-localrep}(iii), it is mutually local with respect to all fields from ${\rm Fld}(V)$. One also has
 $$
X(z)|0\rangle =  \partial_z^{(m)} \left(\sum_{n=0}^\infty T^{(n)} a z^n\right) = \sum_{n=0}^\infty  \binom{m+n}{m} T^{(m+n)} a z^n = \sum_{n=0}^\infty  T^{(n)}(T^{(m)}a)z^n\ .
$$
The latter is $Y(T^{(m)} a,z)|0\rangle$, and so Lemma \ref{lem-Goddard} implies that $X(z)=Y(T^{(m)} a,z)$.
	
(ii) The equality follows from Lemmas \ref{lem-vacuum}(ii) and \ref{lem-fs}.
\end{proof}

\begin{cor}\label{cor-Tderiv}
For all $a,b\in V'$ and $n\in\bfZ$, one has 
$$
T(a_{(n)} b) = (Ta)_{(n)} b + a_{(n)}(Tb)\ .
$$
\end{cor}

\begin{proof}
Applying the operator $fs$ to the equality of Lemma \ref{lem-partialn}(i) and using the equalities of Lemma \ref{lem-basic}(ii), we get the required fact.
\end{proof}

\begin{cor}\label{cor-prodnorm}
For every pair $(a,b)\in V'\times V$ and every $n\in\bfZ$, one has $||a_{(n)}b||\le||a||'\cdot ||b||$. If in addition, $b\in V'$, then $||a_{(n)}b||'\le||a||'\cdot||b||'$.
\end{cor}

\begin{proof}
One has $||a_{(n)}b||\le ||a_{(n)}||\cdot||b||\le ||a||'\cdot ||b||$. If $b\in V'$, then $a_{(n)}b\in V'$ as well, and one has 
 $$
||a_{(n)}b||' = ||Y(a_{(n)}b,z)|| \le ||Y(a,z)||\cdot ||Y(b,z)|| = ||a||'\cdot||b||'\ .
 $$
The statement follows.
\end{proof}

\begin{cor}\label{cor-Borcherds}(Borcherds identity) For all $a,b\in V'$ and $n\in\bfZ$, setting $a(z)=Y(a,z)$ and $b(z)=Y(b,z)$, one has
	$$
	a(z)b(w)i_{z,w}((z-w)^n) - b(w)a(z)i_{w,z}((z-w)^n) = \sum_{j=0}^\infty (a_{(n+j)} b)(w) \partial^{(j)}_w \delta(z-w)\ . 
	$$
\end{cor} 

\begin{proof}
	The identity follows from the equality of Corollary \ref{cor-preBorcherds} in which, by Lemma \ref{lem-basic}(ii), $a(w)_{(n+j)}b(w)$ is replaced by  $(a_{(n+j)} b)(w)$.
\end{proof}

\begin{cor}\label{cor-Liev}
The $K$-submodule of $\End(V)$, generated by the operators $a_{(n)}$ for all $a\in V'$ and $n\in\bfZ$, is a Lie subalgebra of $\End(V)$, denoted by ${\rm Lie}_V$.
\end{cor}

\begin{proof}
The statement follows from Lemmas \ref{lem-localrep}(ii) and \ref{lem-basic}(ii).	
\end{proof}

\begin{lem}\label{lem-skew}
For all $a,b\in V'$, one has:
\begin{itemize}
    \item [(i)]
$e^{wT}	Y(a,z) e^{-wT} = i_{z,w} Y(a,z+w)$;
    \item[(ii)] (Skewsymmetry)
$Y(a,z) b = e^{zT} Y(b,-z) a$.
\end{itemize}	
\end{lem}

\begin{proof}
(i) Let $f(z,w)$ and $g(z,w)$ denote the left and right hand sides of the equality. One has $f(z,0) = Y(a,z)=g(z,0)$ and, therefore, 
 $$
\partial_w f(z,w) = T e^{wT} Y(a,z) e^{-wT} - e^{wT} Y(a,z) T e^{-wT} = [T,f(z,w)]. 
 $$
On the other hand, by the translation covariance axiom, one has $\partial_w g(z,w) =[T,g(z,w)]$, and the required fact follows from Lemma \ref{lem-diffequ}(i).
 \medskip

(ii) Mutual locality of $Y(a,z)$ and $Y(b,z)$ implies that
$$
(z-w)^N Y(a,z)e^{wT} b - (z-w)^N Y(b,w)e^{zT} a \to 0\ {\rm as\ } N\to+\infty\ .
$$
By (i), the second summand is equal to 
 $$
(z-w)^N e^{zT} e^{-zT} Y(b,w) e^{zT} a =  e^{zT} i_{w,z} ((z-w)^N Y(b,w-z)) a\ . 
 $$
Since $Y(b,w-z) \in \End(V)\langle\langle w-z \rangle\rangle$, it follows that
 $$
(z-w)^N Y(b,w-z) - i_{w,z} ((z-w)^N Y(b,w-z))\to 0\ {\rm as\ }N\to +\infty
 $$
and, therefore, one has
 $$
(z-w)^N Y(a,z)e^{wT} b - (z-w)^N e^{zT} Y(b,w-z)a \to 0\ {\rm as\ }N\to +\infty\ .   
 $$
Substituting $w=0$ and using the fact that multiplication by $z^N$ is an isometric operator, we get the required fact.
\end{proof}

\begin{cor}\label{cor-skew}
(i) For all $a,b\in V'$ and $n\in\bfZ$, one has
 $$
a_{(n)} b = -\sum_{m=0}^\infty (-1)^{m+n} T^{(m)} (b_{(m+n)} a)\ ; 
 $$	
 
(ii) for all $a(z),b(z)\in{\rm Fld}(V)$, one has
 $$
a(z)_{(n)}b(z) = - \sum_{m=0}^\infty (-1)^{m+n} \partial^{(m)}_z (b(z)_{(m+n)} a(z))\ . 
 $$
\end{cor}

\begin{proof}
 The statement (i) straightforwardly follows from Lemma \ref{lem-skew}(ii). The statement (ii) follows from (i) and Lemma \ref{lem-basic}(ii).
\end{proof}

In the same way Corollary \ref{cor-Borcherds} implies two other forms of the Borcherds identity.

\begin{cor}
  \label{cor-B}
(i) For all $a,b,c\in V'$ and $m,n,k\in\bfZ$, one has
$$
  \sum_{j=0}^\infty \binom{m}{j}(a_{(n+j)}b)_{(m+k-j)}c =
  \sum_{j=0}^\infty (-1)^j \binom{n}{j} ((a_{(m+n-j)}(b_{(k+j)}c)
  -(-1)^n(b_{(n+k-j)}(a_{(m+j)}c));
$$

  (ii) identity from (i) holds if we replace all $a,b,c\in V'$  by
  $a(z),b(z),c(z)\in{\rm Fld}(V)$.    
\end{cor}

Note that Corollary \ref{cor-skew} follows from Corollary \ref{cor-B}
by letting $c=\ket{0}$, $m=-1$ and $k=0$. 

\subsection{Commutative vertex $K$-algebras}

\begin{lem}\label{lem-commut}
	The following properties of a vertex $K$-algebra $V$ are equivalent:
	\begin{itemize}
		\item [(a)]
		$\varphi(z)\in \End(V)[[z]]$ for all $\varphi(z)\in{\rm Fld}(V)$;
		\item [(b)]
		$\varphi(z)\in \End(V)((z))$ for all $\varphi(z)\in{\rm Fld}(V)$;
		\item[(c)]
		$[\varphi(z),\psi(w)]=0$ for all $\varphi(z),\psi(z)\in{\rm Fld}(V)$.	
	\end{itemize}
\end{lem}

\begin{proof}
The implication (a)$\Longrightarrow$(b) is trivial. Suppose (b) is true. Then for all $\varphi(z),\psi(z)\in{\rm Fld}(V)$, one has $[\varphi(z),\psi(w)] \in \End(V)[[w^{\pm 1}]]((z))$. Since the fields $\varphi(z)$ and $\psi(z)$ are mutually local, $(z-w)^N [\varphi(z),\psi(w)]\to0$ as $N\to\infty$. But the multiplication by $z-w$ on $\End(V)[[w^{\pm 1}]]((z))$ is an isometry, by Lemma \ref{lem-zminw}. This implies (c). Finally, assume (c) is true. Then $\varphi(z)\psi(w)|0\rangle\in V[[z^{\pm 1},w]]$ and $\psi(w)\varphi(z)|0\rangle\in V[[z,w^{\pm 1}]]$. Since both vectors are equal, it follows that they both are contained in $V[[z,w]]$. Substituting $w=0$, we get $\varphi(z)\psi_{(-1)}|0\rangle\in V[[z]]$ for all $\psi(z)\in{\rm Fld}(V)$. Since the $K$-saturation of the image $V'$ of the homomorphism $fs$ is dense in $V$, we get $\varphi(z)\in \End(V)[[z]]$.
\end{proof}

A vertex $K$-algebra $V$ that possesses the equivalent properties of
Lemma \ref{lem-commut} is called {\it commutative}. For example, the vertex $K$-algebras with zero translation operator (from Example \ref{exam-T0}) are commutative.

\begin{defin}
A {\it strongly commutative Banach $K$-algebra} is a quadruple $(V,|0\rangle,T,A)$ consisting of a Banach torsion free $K$-module $V$, an element $|0\rangle\in V$, a bounded operator $T\in\End(V)$ with $T|0\rangle=0$, and a Banach $K$-subalgebra $A\subset \End(V)[[z]]$ such that 
\begin{itemize}
	\item [(C.1)]
$[a(z),b(w)]=0$ for all $a(z),b(z)\in A$;
    \item[(C.2)]
$A$ is preserved by the operators $\partial^{(n)}_z$ for all $n\ge0$;
    \item[(C.3)]
$[T,a(z)] = \partial_z a(z)$ for all $a(z)\in A$;
    \item[(C.4)]
if $V'$ is the image of the homomorphism $A\to V:a(z)\mapsto a(0)|0\rangle$, then the $K$-saturation of $V'$ is dense in $V$.
\end{itemize}
\end{defin}

Notice that the homomorphism in (C.4) is injective (cf. Lemma \ref{lem-fs}). Indeed, suppose that for $a(z)=\sum_{n=0}^\infty a^{(n)} z^n\in A$, one has $a(0)|0\rangle=a^{(0)}|0\rangle=0$. The property (C.3) implies that $T a^{(n)}|0\rangle = (n+1) a^{(n+1)}|0\rangle$ for all $n\ge0$ and, by induction, we get $a^{(n)}|0\rangle=0$ for all $n\ge0$, i.e., $a(z)|0\rangle=0$. Then for every $b(z)\in A$, one has $a(z)b(w)|0\rangle = b(w)a(z)|0\rangle=0$. It follows that $a(z)b(0)|0\rangle=0$ for all $b(z)\in A$, i.e., $a(z)V'=0$, and the property (C.4) implies that $a(z)=0$. Thus, we can denote by $a(z)$ the preimage of an element $a\in V'$ in $A$.

A morphism of such objects $(V,|0\rangle,T,A) \to (U,|0\rangle,T,B)$ is a bounded $K$-linear homomorphism $\varphi:V\to U$ which takes $|0\rangle$ to $|0\rangle$ and $V'$ to $U'$, commutes with $T$, and the induced map $V'\to B: a\mapsto \varphi(a)(z)$ gives rise to a bounded homomorphism of Banach $K$-algebras $A\to B$.

\begin{lem}
The correspondence $V\mapsto (V,|0\rangle,T,{\rm Fld}(V))$ gives rise to an equivalence between the category of commutative vertex $K$-algebras and the category of strongly commutative Banach $K$-algebras.
\end{lem}

\begin{proof}
First of all, we notice that for any fields $\varphi(z),\psi(z) \in \End(V)[[z]]$, one has $\varphi(z)_{(n)}\psi(z)=0$ for $n\ge0$, and $\varphi(z)_{(n)}\psi(z)= (\partial^{(m)}_z \varphi(z))\psi(z)$ for $n=-m-1<0$. Let $V$ be commutative vertex $K$-algebra. Then the property (C.1) holds, by Lemma \ref{lem-commut}, and (C.3) and (C.4) hold by (V.3) and (V.5), respectively. The above remark implies that ${\rm Fld}(V)$ is preserved by multiplication and by the operators $\partial^{(n)}_z$, i.e., $(V,|0\rangle,T,{\rm Fld}(V))$ is a strongly commutative Banach $K$-algebra. Conversely, if the latter holds, then all of the properties (V.1)-(V.5) evidently hold, and so $V$ is a commutative vertex $K$-algebra. That the correspondence considered induces an equivalence of categories follow easily from the definitions of morphisms in both categories.
\end{proof}

\begin{cor}
Let $V$ be a commutative vertex $K$-algebra. Then
\begin{itemize}
	\item [(i)]
the $K$-bilinear operation $V'\times V\to V: (a,b)\mapsto a\cdot b=Y(a,0)b$ defines the
structure of a commutative  Banach $K$-algebra on $(V',||\ ||')$ with the unit
element $|0\rangle$, and the structure of a Banach module on $V$ over the latter;
    \item[(ii)]
for every pair $(a,b)\in V'\times V$, one has $Y(a,z)b = (e^{zT} a)\cdot b$;
   \item[(iii)]
the operator $T$ on $(V',||\ ||')$ is a derivation with respect to the above multiplication and of norm at most one.
\end{itemize}	
\end{cor}

In the admissible case, $V'$ is a Banach $K$-algebra with respect to the norm $||\ ||$ since it is equivalent to the norm $||\ ||'$.

\begin{proof}
(i) Let $a,b\in V'$. Since $b=Y(b,0)|0\rangle$, we have $a\cdot b=Y(a,0)Y(b,0)|0\rangle$ and, by the property (c) of Lemma \ref{lem-commut}, $a\cdot b=b\cdot a$. Furthermore, the equality $Y(a,0)Y(b,0) = Y(b,0)Y(a,0)$ gives the
equality $a\cdot (b\cdot c) = b\cdot(a\cdot c)$ for all $a,b,c\in
V'$, and we get $a\cdot (b\cdot c) = a\cdot (c\cdot b) = c\cdot
(a\cdot b) = (a\cdot b)\cdot c$. If $a\in V'$ and $b\in V$, then $||a\cdot b||\le
||a||'\cdot||b||$. If in addition $b\in V'$, then $||a\cdot b||' = ||a_{(-1)}b||'\le ||a||'\cdot||b||'$, by Corollary \ref{cor-prodnorm}.
 \medskip 

(ii) Since $Y(a,z) = \sum_{n=0}^\infty T^{(n)} a z^n$, one has
 $$
Y(a,z)b = \sum_{n=0}^\infty (T^{(n)} a) b z^n = \sum_{n=0}^\infty Y(T^{(n)} a,0)b z^n = \sum_{n=0}^\infty ((T^{(n)} a)\cdot b) z^n = (e^{zT}a) \cdot b\ ,
 $$
and the statement (ii) follows.
 \medskip

(iii) For $a,b\in V'$, Lemma \ref{cor-Tderiv} implies that
 $$
T(a\cdot b) = T(a_{(-1)}b) = (Ta)_{(-1)} b+a_{(-1)}(Tb) =  (Ta)\cdot b+a\cdot(Tb)\ ,
 $$
i.e., $T$ is a derivation. Using (i), we get $||Ta||'=||a_{(-2)}|0\rangle||' \le ||a||'\cdot |||0\rangle||'\le||a||'$ (if $|0\rangle\not=0$, the second inequality is an equlity). Thus, the operator $T$ on $(V',||\ ||')$ is of norm at most one.
\end{proof}

\begin{exams}\label{exam-vertex}
(i) Let $A$ be a commutative Banach $K$-algebra provided with a system of bounded $K$-linear operators $\{T^{(n)}\}_{n\ge0}$ such that $T^{(0)}=I_A$ and $T^{(n)}(a\cdot b) = \sum_{i+j=n} T^{(i)} a \cdot T^{(j)} b$ for all $n\ge0$ and $a,b\in A$. Assume that the set $A'=\{a\in A \big\vert$ there exists $C>0$ with  $||T^{(n)} a||\le C||a||$ for all $n\ge0\}$ is dense in $A$. Then $A$ is a commutative vertex $K$-algebra with $Y(a,z)$, defined for $a\in A'$ as multiplication by $e^{zT} a = \sum_{n=0}^\infty T^{(n)} a z^n$. If there exists $C>0$ with $||T^{(n)}||\le C$ for all $n\ge0$, then $A$ is admissible.

(ii) The commutative Banach $K$-algebra $K\{t^{\pm 1}\}$, provided with the system of operators $T^{(n)}=\partial^{(n)}_t$, $n\ge0$, is an admissible commutative vertex $K$-algebra.

(iii) For a positive real number $r$, let $K\{r^{-1} t\}$ denote the commutative Banach $K$-algebra of formal power series $f=\sum_{n=0}^\infty \lambda_n t^n$ with $|\lambda_n| r^n\to0$ as $n\to\infty$, provided with the Banach norm $||f||=\max_n\{|\lambda_n| r^n\}$. The derivation operator $\partial_t$ has the norm $r^{-1}$ and, for $n\ge0$, the operator $T^{(n)}=\partial_t^{(n)}$ has the norm $r^{-n}$. Thus, if $r\ge1$, then $K\{r^{-1} t\}$ is an admissible commutative vertex $K$-algebra. If $r<1$, $A$ is not admissible (see Example \ref{exam-nonadmis}).
\end{exams}

\subsection{Extension Theorem} \label{subsec-extheor}
Suppose we are given a torsion free Banach $K$-module $V$, an element $|0\rangle\in V$, a bounded endomorphism $T\in\End(V)$, and a collection of fields
 $$
\calF = \left\{ a^j(z) = \sum_{n\in\bfZ} a^j_{(n)} z^{-n-1}\right\}_{j\in J} \subset  {\rm End}(V)\langle\langle z\rangle\rangle\ .
 $$ 
Let $V_\calF$ be the minimal $K$-submodule of $V$ that contains the element $|0\rangle$ and is preserved under $K$-linear endomorphisms of the form $v\mapsto a^j_{(n)} v$ for $j\in J$ and $n\in\bfZ$. Furthermore, suppose that the above data possess the following properties:
\begin{itemize}
	\item [(E.1)]
	(vacuum axiom) $T|0\rangle=0$ and $a^j(z)|0\rangle \in V[[z]]$ for all $j\in J$;
	\item[(E.2)]
	(translation covariance) $[T,a^j(z)]=\partial_z a^j(z)$ for all $j\in J$;
	\item[(E.3)]
	(locality) the fields $a^j(z)\in\calF$ are mutually local;
	\item[(E.4)]
	(completeness) the $K$-saturation of the $K$-submodule $V_\calF$ is dense in $V$.
\end{itemize}

\begin{theor}\label{theor-gener}
In the above situation, there is a vertex $K$-algebra structure on $V$ in which ${\rm Fld}(V)$ is the closed $K$-submodule of ${\rm End}(V)\langle\langle z\rangle\rangle$, generated by $\calF$.
\end{theor}

\begin{proof}
Validity of the properties (V.1)-(V.4) from Definition \ref{def-vertex} follows from Lemma \ref{lem-localclosed}, and validity of the property (V.5) from the same definition follows from Lemma \ref{lem-vacuum}(ii).
\end{proof}

Let $U$ be a vertex $K$-algebra, and let $V$ be a vertex $K'$-algebra $V$ for a commutative non-Archimedean Banach $K$-algebra $K'$. The completed tensor product $U\hatotimes_K V$ is a Banach $K'$-module. Assume that $K'$ and $U\hatotimes_K V$ possess the property from the beginning of this section, i.e., 
they are torsion free abelian groups, and assume that $U\hatotimes_K V$ is a torsion free $K'$-module.

Notice that the valuation on $\bfZ$ for $K'$ is not necessarily the same as for $K$. For example, they do not coincide for $K=\bfZ$ with the trivial valuation and $K'=\bfZ_p$, the ring of $p$-adic integers, with a $p$-adic valuation.

\begin{cor}\label{cor-tensor}
In the above situation, the completed tensor product $U\hatotimes_K V$ is a vertex $K'$-algebra with the vacuum vector $|0\rangle\hatotimes|0\rangle$, the translation operator $T\hatotimes I_{V} + I_U \hatotimes T$, and the space of states ${\rm Fld}(U\hatotimes_K V)$, which is the closed $K$-submodule of $\End(U\hatotimes_K V) \langle\langle z\rangle\rangle$, generated by fields of the form $\varphi(z)\hatotimes\psi(z)$ for $\varphi(z)\in{\rm Fld}(U)$ and $\psi(z)\in {\rm Fld}(V)$.
 \hfill$\square$
\end{cor}

\begin{exam}\label{exam-nonadmis}
Assume that $r<1$ in Example \ref{exam-vertex}. Then the norm of $\partial_t^{(n)}$ tends to infinity together with $n$. We apply Theorem \ref{theor-gener} to the Banach space $V=K\{r^{-1} t\}$, the vacuum element $1\in K\{r^{-1} t\}$, the translation operator $T=\partial_t$, and the collection of fields $\calF=\{e^{zT} a \big\vert a\in K[t]\}$. The properties (E.1)-(E.3) clearly hold, and the property (E.4) holds because the ring of polynomials $K[t]$ is dense in $K\{r^{-1} t\}$. Thus, $K\{r^{-1} t\}$ is a commutative vertex $K$-algebra, which is not admissible.
\end{exam}

In all of the examples of vertex $K$-algebras, considered in the next subsections, we start from an ordinary vertex algebra $V$ over a commutative ring $R$, which is torsion free as an abelian group, with both $R$ and $V$ provided with the trivial Banach norm. Then for each injective homomorphism $R\to K$ to a commutative non-Archimedean ring $K$, by Corollary \ref{cor-tensor}, we get a non-Archimedean vertex $K$-algebra $V_K=V\hatotimes_R K$. 

\begin{lem}\label{lem-fsurj}
In the above situation, assume that $V$ is a free $R$-module. Then the following is true:
\begin{itemize}
\item[(i)]
if the field-state homomorphism $fs:{\rm Fld}(V)\to V$ is bijective, the field-state homomorphism ${\rm Fld}(V_K)\to V_K$ is an isometric isomorphism;
\item[(ii)]
if the Banach norm on $K$ induces the trivial norm on $R$, then the field-state homomorphism ${\rm Fld}(V_K)\to V_K$ is an isometric injection.
\end{itemize}
In particular, in both cases the vertex $K$-algebra $V_K$ is admissible.
\end{lem}

\begin{proof}
Let $\{e_i\}_{i\in I}$ is a basis of $V$ over $R$. Each element $v\in V_K$ has a unique representation as a sum $v=\sum_{i\in I} \lambda_i e_i$ with $\lambda_i\in K$ and $\lambda_i\to0$ with respect to the filter of complements of finite subsets of $I$, and one has $||v||=\max_{i\in I} \{||\lambda_i||\}$. This implies that the canonical homomorphisms $V\to V_K$ and $\End(V)\to\End(V_K)$ are isometric. Furthermore, Let $U$ be the $K$-submodule of $V$ consisting of finite sums $v=\sum_{i\in I} \lambda_i e_i$. Then $U$ is dense in $V_K$ in the case of (i), and its $\bfZ$-saturation is dense in $V_K$ in the case of (ii). 

(i) For an element $v=\sum_{i\in I} \lambda_i e_i\in U$, we set $\varphi(z)=\sum_{i\in I} \lambda_i  sf(e_i)\in\End(V_K)((z))$. Since $fs(\varphi(z))=v$, one has $||f||\le||\varphi(z)||$. On the other hand, since $||sf(e_i)||=1$, one has $||\varphi(z)||\le \max_{i\in I} \{||\lambda_i||\}=||||v||$. This implies that $||fs(\varphi(z))||=||\varphi(z)||$, and the required statement follows.

(ii) For each nonzero  element $v=\sum_{i\in I} \lambda_i e_i\in U$ there exists $N\ge1$ with $Ne_i\in V'$ for all $i\in I$ with $\lambda_i\not=0$ and, in particular, $Nv\in V'$. We set $\varphi(z)=\sum_{i\in I} \lambda_i sf(Ne_i)$. Then $\varphi(z)\in{\rm Fld}(V_K)$ and $fs(\varphi(z)) = N\cdot \sum_{i\in I} \lambda_i e_i$.  By the assumption, $||Ne_i||=1$ and, therefore, one has
 $$
||\varphi(z)|| \le \max_{i\in I} \{||\lambda_i||\}= ||fs(\varphi(z))|| \le||\varphi(z)||\ .
 $$
Thus,  $||fs(\varphi(z))|| = ||\varphi(z)||$. This implies the required fact.
\end{proof}

\subsection{Vertex $K$-algebra associated to Lie algebras}\label{subsec-lieassoc}
Let $R$ be a commutative ring, torsion free as an abelian group, let $\gothg$ be a Lie $R$-algebra, torsion free as an $R$-module, let $T$ be a derivation of $\gothg$, and let $\calF$ be a collection of mutually local {\it formal distributions} 
 $$
\calF = \left\{ a^j(z) = \sum_{n\in\bfZ} a^j_{(n)} z^{-n-1}\right\}_{j\in J} \subset \gothg[[z^{\pm1}]]\ .
 $$
The action of $T$ on $\gothg$ extends naturally to an action on the space of formal distributions $\gothg[[z^{\pm1}]]$. Let $\calF_T$ denote the minimal $R$-submodule of $\gothg[[z^{\pm1}]]$ that contains $\calF$ and is invariant under the action of $T$. Suppose also that the above data possess the following properties:
\begin{itemize}
	\item [(L.1)]
the $R$-saturation of the $R$-submodule of $\gothg$, generated by the elements $a^j_{(n)}$ for all $j\in J$ and $n\in\bfZ$, coincides with $\gothg$;
	\item [(L.2)]
$T(a^j(z)) = \partial_z a^j(z)$ for all $j\in J$;
    \item [(L.3)]
for every $n\ge0$, the $n$-th product of any pair of formal distributions from $\calF$ lies in $\calF_T$.
\end{itemize}

Let $\gothg_{-}$ be the $R$-submodule of $\gothg$, generated by the elements $a^j_{(n)}$ for all $j\in J$ and $n\ge0$. It follows from the property (L.2) and Lemma \ref{lem-localrep}(ii)  that $\gothg_{-}$ is a $T$-invariant Lie $R$-subalgebra of $\gothg$ (the {\it annihilation subalgebra}). Let $V$ be the quotient $U(\gothg)/\bfa$ of the universal enveloping $R$-algebra $U(\gothg)$ of $\gothg$ by the $R$-saturation $\bfa$ of the left ideal $U(\gothg)\gothg_{-}$. 

Notice that left miltiplication by elements of $\gothg$ in $U(\gothg)$ induces an $R$-linear action of $\gothg$ on $V$. This gives rise to a homomorphism $\gothg[[z^{\pm1}]]\to \End[[z^{\pm1}]]$, and it takes the collection $\calF$ to a collection of mutually local formal distributions $\calF\big\vert_V$ in $\End[[z^{\pm1}]]$. Notice also that the derivation $T$ on $\gothg$ induces an action of $T$ on $V$. 

\begin{lem}\label{lem-Lievertex}
In the above situation, formal distributions in $\calF\big\vert_V$ are in fact $\End(V)$-valued quantum fields, i.e., $\calF\big\vert_V \subset \End(V)\langle\langle z \rangle\rangle$. 	
\end{lem}

Notice that since the valuation on $V$ is trivial, then $\End(V)\langle\langle z \rangle\rangle = \End(V)((z))$.

\begin{proof}
Step 1. Let $F$ be the $R$-submodule of formal distributions $a(z)\in \gothg[[z^{\pm1}]]$ with $a_{(n)}\in \gothg_{-}$ for all $n\ge0$. By the definition, one has $\calF\subset F$, and the property (L.2) implies that $\calF_T \subset F$. From the property (L.3) it then follows that, for every $i\ge0$ and every pair of formal distributions $a(z),b(z)\in\calF$, one has $a(z)_{(i)}b(z) \in F$ and, in particular, $(a(z)_{(i)}b(z))_{(k)}\in\gothg_{-}$ for all $k\ge0$.
 \medskip

Step 2. {\it For any pair $a(z),b(z) \in \calF$ and $n\in\bfZ$, one has} 
 $$ 
[a(z),b_{(n)}] \in \gothg((z))+\gothg_{-}[[z^{\pm1}]] \subset \gothg[[z^{\pm1}]]\ .
 $$ 
Indeed, by Lemma \ref{lem-localrep}(ii), one has
 $$
[a(z),b_{(n)}] = \sum_{m\in\bfZ} [a_{(m)},b_{(n)}] z^{-m-1} =  \sum_{m\in\bfZ} \sum_{i=0}^\infty \binom{m}{i} (a(z)_{(i)} b(z))_{(m+n-i)} z^{-m-1}\ .
 $$
Since $a(z)$ and $b(z)$ are mutually local, there exists $N\ge0$ with $a(z)_{(i)} b(z)=0$ for all $i>N$. Thus, for $m\ge0$ the coefficient at $z^{-m-1}$ is equal to
 $$
\sum_{i=0}^{\min({m,N)}} \binom{m}{i} (a(z)_{(i)} b(z))_{(m+n-i)}\ .
 $$
If $n\ge0$, then $m+n-i\ge0$ for all $0\le i\le m$ and, by Step 1, that coefficient belongs to $\gothg_{-}$. If $n=-k-1$ for $k\ge0$, then for $m\ge N+k+1$, one has $m+n-i\ge0$ and, by Step 1, the coefficient at $z^{-m-1}$ belongs to $\gothg_{-}$.
 \medskip
 
Step 3. {\it The statement of the lemma is true.} Indeed, by the property (L.1), it suffices to show that, for any $b(z)\in\calF$ and any element $x\in U(\gothg)$ of the form $x=a^{j_1}_{(n_1)} \cdot\ldots\cdot a^{(j_s)}_{(n_s)}$, one has $b(z)x\in U(\gothg)((z)) + U(\gothg)g_{-}[[z^{\pm1}]]$. If $s=0$, then $x=1$ and the required fact follows from the inclusion $a(z)\in F$ (see Step 1). Assume that $s\ge1$ and the required fact is true for smaller values of $s$. Then $x=a^j_{(n)}y$ for $y$ of the same form with smaller value of $s$. One has
 $$
b(z)x = [b(z),a^j_{(n)}]y +  a^j_{(n)} b(z) y\ .
 $$
By Step 2, the first summand on the right hand side belongs to $U(\gothg)((z)) + U(\gothg)g_{-}[[z^{\pm1}]]$ and, by the induction, the same holds for the second summand. 
\end{proof}

\begin{cor}\label{cor-Lievertex}
In the above situation, there is a vertex $R$-algebra structure on $V$ in which the vacuum vector $|0\rangle$ is the image of $1$ in $V$, the endomorphism $T$ of $V$ is induced by the derivation $T$ on $\gothg$, and ${\rm Fld}(V)$ is the $R$-submodule of $\End(V)\langle\langle z\rangle\rangle$, generated by $\calF\big\vert_V$.
 \hfill$\square$
\end{cor}

The vertex $R$-algebra from Corollary \ref{cor-Lievertex} is denoted by $V(\gothg,\calF)$.

\subsection{Examples of vertex $K$-algebras}\label{subsec-bosfer}
\ \medskip

\centerline{\it Free boson vertex $K$-algebra}
\medskip
Let $B=\bfZ[x_1,x_2,\ldots]$ be the ring of polynomials over $\bfZ$ in the variables $x_1,x_2,\ldots$ with $\bfZ$ and $B$ provided with the trivial Banach norm. Let also $|0\rangle=1$, $T=\sum_{i=2}^\infty (i-1)x_i {\partial\over{\partial x_{i-1}}}: B\to B$, and $\calF=\{a(z)\}$ for $a(z)= \sum_{n\in\bfZ} a_{(n)} z^{-n-1}$, defined by $a_{(n)} =n {\partial\over{\partial x_n}}$ for $n>0$, $a_{(n)} = x_{-n}$ for $n<0$, and $a_{(0)}=0$. These data possess the properties (E.1)-(E.4) from \S \ref{subsec-extheor} and so, by Theorem \ref{theor-gener}, they give rise to an admissible vertex $\bfZ$-algebra structure on $B$. Notice that $B$ is a free abelian group that consists of {\it finite} sums $\sum_\mu \lambda_\mu x^\mu$, taken over sequences of non-negative integers $\mu=(\mu_1,\mu_2,\ldots)$ all but finitely many of which are equal to zero, and where $x^\mu=\prod_{i=1}^\infty x^{\mu_i}$ and $\lambda_\mu\in \bfZ$. Notice also that, since the $K$-submodule $B_\calF$ from (E.4) coincides with $B$, the field-state homomorphism $fs: {\rm Fld}(B)\to B$ is bijective. Thus, by Lemma \ref{lem-fsurj}(i), for any $K$, $B_K=B\hatotimes_\bfZ K$ is an admissible vertex $K$-algebra for which the field-state homomorphism $fs:{\rm Fld}(B_K)\to B_K$ is an isometric isomorphism.

Sometimes in literature, the free boson vertex algebra is presented in a slightly different form. Namely, let $B^t=\bfZ[y_1,y_2.\ldots]$ be a similar ring of polynomials, and let $|0\rangle=1$, $T=\sum_{i=2}^\infty iy_i {\partial\over{\partial y_{i-1}}}: B\to B$, and $\calF^t=\{b(z)\}$ for $b(z)= \sum_{n\in\bfZ} b_{(n)} z^{-n-1}$, defined by $b_{(n)} ={\partial\over{\partial y_n}}$ for $n>0$, $b_{(n)} = -n y_{-n}$ for $n<0$, and $b_{(0)}=0$. Notice that these data possess the same properties (E.1)-(E.4) from \S \ref{subsec-extheor} and so, by Theorem \ref{theor-gener}, they give rise to an admissible vertex $\bfZ$-algebra structure on $B^t$. But in this case, the field-state homomorphism $fs:{\rm Fld}(B^t) \to B^t$ is not surjective. There is a homomorphism of vertex $\bfZ$-algebras $B\to B^t: x_i\mapsto iy_i$, which becomes an isomorphism after tensoring with $\bfQ$. 

In any case, by Corollary \ref{cor-tensor}, for any $K$, $B^t_K$ is a vertex $K$-algebra. By Lemma \ref{lem-fsurj}(ii), if the valuation on the image of $\bfZ$ in $K$ is trivial, then $B^t_K$ is admissible. Assume now that the valuation $|\ |$ on the image of $\bfZ$ in $K$ is $p$-adic. {\it We claim that in this case $B^t_K$ is not admissible.} Indeed, one has $||b_{(n)}||=1$ for $n>0$ and $||b_{(n)}||=|n|$ for $n<0$. The space ${\rm Fld}(B^t_K)$ contains also the fields $\partial^{(m)}_z b(z)$ for all $m\ge1$. One has
$$
\partial^{(m)}_z b(z) = \sum_{n\in\bfZ} (-1)^m\binom{m+n}{m} b_{(n)} z^{-m-n-1} 
$$
and, therefore, $fs(\partial^{(m)} b(z)) = b_{(-m-1)}|0\rangle = (m+1) y_{m+1}$. The norm of the latter element is equal to $|m+1|$. Consider the above element for $m=p^k-1$ with $k\ge1$. For this value of $m$, one has $||fs(\partial^{(m)}_z b(z))|| = |p|^k$. On the other hand, the binomial coefficient in the above sum for $n=p$ is equal to $\binom{m+p}{p}$ and, therefore, it is not divisible by $p$. Since $||b_{(p)}||=1$, this implies that $||\partial^{(m)}_z b(z)||=1$. Since $|p|^k\to0$ as $k\to\infty$, we see that the operator $fs$ is not admissible and, therefore, the vertex $K$-algebra $B^t_K$ is not admissible.
 \medskip

\centerline{\it Free fermion vertex $K$-algebra}
\medskip

For simplicity of exposition, we were considering the spaces of states with the trivial parity $p=\overline0$. This discussion easily extends to the super spaces. Namely, $V=V_{\overline0}\oplus V_{\overline1}$, where $V_{\oj}$ for $\oj\in\bfZ/2\bfZ$ are $K$-modules, and the parity $p$ is defined by $p(v)=\oj$ for $v\in V_{\oj}$. In the vertex algebra definition, one assumes that $|0\rangle\in V_{\overline0}$, $T V_\oj\subset V_\oj$, and the commutator is understood in the ``super" sense. The free fermion algebra is the most important example of such an object.

Let $F=\bfZ[\xi_1,\xi_2,\ldots]$ be the Grassmann superalgebra over $\bfZ$ in $\xi_1,\xi_2,\ldots$ with $\xi_i\xi_j=-\xi_j\xi_i$, $\xi^2_i=0$, and $p(\xi_i)=\overline 1$. Both $\bfZ$ and $F$ are provided with the trivial norm. Let also $|0\rangle=1$ and $T=\sum_{i=1}^\infty i\xi_{i+1} {\partial\over{\partial \xi_i}}$, where $\partial\over{\partial \xi_i}$ is the odd derivation (i.e., ${{\partial(fg)} \over {\partial \xi_i}} = {{\partial f} \over {\partial \xi_i}} g + (-1)^{p(f)} f {{\partial g}\over{\partial\xi_j}}$) such that ${{\partial \xi_i}\over{\partial\xi_j}} = \delta_{i,j}$. Finally, let $\calF=\{\varphi(z)\}$ for $\varphi(z)= \sum_{n\in\bfZ} \varphi_{(n)} z^{-n-1}$, defined by $\varphi_{(n)} = {\partial\over {\partial\xi_{n+1}}}$ for $n\ge0$, and $\varphi_{(n)} = \xi_{-n}$ for $n<0$. These data data possess the properties (E.1)-(E.4) from \S \ref{subsec-extheor} and so, by Theorem \ref{theor-gener}, they gives rise to an admissible vertex $\bfZ$-algebra structure on $F$. Notice that $F$ is the free abelian group that consists of {\it finite} sums $\sum_\nu a_\nu \xi^\nu$, taken over finite subsets $\nu=\{\nu_1<\ldots<\nu_m\}\subset \bfZ_{>0}$, and where $\alpha_\nu\in\bfZ$, $\xi^\nu = \xi_{\nu_1} \wedge\ldots\wedge \xi_{\nu_m}$, and $\xi^\emptyset=1$. Notice also that, in this case, the field-state homomorphism $fs$ is a bijection. Thus, by Lemma \ref{lem-fsurj}(i), for any $K$, $F_K$ is an admissible vertex $K$-algebra for which the field-state homomorphism $fs:{\rm Fld}(F_K)\to F_K$ is an isometric isomorphism.
 \medskip
 
\centerline{\it Universal Virasoro vertex $K$-algebra}
 \medskip
 
The Virasoro algebra is a Lie algebra Vir over the ring $R=\bfZ[\frac{1}{2}]$, the localization of $\bfZ$ by powers of $2$, with free generators $L_n$ for $n\in\bfZ$ and $C$ satisfying the following relations
 $$
[L_m,L_n] = (m-n) L_{m+n} + \delta_{m,-n} \frac{m^3-m}{12} C,\ [C,L_m]=0\ {\rm for\ all\ } m,n\in\bfZ\ . 
 $$
(Both $R$ and Vir are provided with the trivial norm.) Setting $L_{(n)}=L_{n-1}$, consider the formal distribution $L(z) = \sum_{n\in\bfZ} L_{(n)} z^{-n-1}$. The above relations can be written in the equivalent form
 $$
[L(z),L(w)] = \partial_w L(w)\delta(z-w) + 2L(w) \partial_w\delta(z-w)  + \frac{C}{2} \partial^{(3)} \delta(z-w)\ .
 $$
It follows that $(z-w)^4 [L(z),L(w)]=0$ and, therefore, $L(z)$ is mutually local with itself. The property (L.1) for $\calF=\{L(z),C\}$ clearly holds. If $T={\rm ad}(L_{-1})$, then $T(L_{(n)}) = [L_{-1},L_{n-1}]= -n L_{(n-1)}$ and, therefore, $T(L(z)) = \partial_z L(z)$, i.e., the property (L.2) holds. The above equality implies that $L(z)_{(0)} L(z) = \partial_z L(z)$, $L(z)_{(1)} L(z) = 2L(z)$, $L(z)_{(3)} L(z) = \frac{C}{2}$, and $L(z)_{(i)} L(z)=0$ for all other values of $i$. This implies that the property (L.3) also holds. Thus, by Corollary \ref{cor-Lievertex}, we get a vertex $R$-algebra ${\rm VIR}=V({\rm Vir},\{L(z),C\})$ called the {\it universal Virasoro vertex algebra}. 

Notice that the annihilation subalgebra of Vir is ${\rm Vir}_{-} = \oplus_{n\ge-1} R L_n$. By the commutation relations, it follows that ${\rm VIR}$ is a free $R$-module that consists of {\it finite} sums $\sum_{m,\mu} \lambda_{m,\mu} C^m L^\mu$, taken over sequences of non-negative integers $(m,\mu)=(m,\mu_1,\mu_2,\ldots)$ all but finitely many of which are equal to zero, and where $L^\mu=\prod_{i=1}^\infty L_{-i-1}^{\mu_i}$ and $\lambda_{m,\mu}\in R$. Notice also that the field-state homomorphism $fs:{\rm Fld}({\rm VIR})\to{\rm VIR}$ is a bijection. Thus, by Lemma \ref{lem-fsurj}, for any $K$ in which $2$ is invertible, ${\rm VIR}_K$ is an admissible vertex $K$-algebra for which the field-state homomorphism ${\rm Fld}({\rm VIR}_K) \to {\rm VIR}_K$ is an isometric isomorphism. Notice that 
 $$
{\rm VIR}_K=\{f= \sum \lambda_{m,\mu} C^m L^\mu \big\vert \lambda_{m,\mu} \in K\ {\rm and\ } ||\lambda_{m,\mu}||\to0\ {\rm as\ } m+\langle\mu\rangle\to\infty\}\ ,
 $$
where the sum is taken over $(m,\mu)$ as above and $\langle\mu\rangle = \sum_{i=0}^\infty i\mu_i$. The Banach norm on ${\rm VIR}_K$ is defined by $||f|| = \max_\mu \{|\lambda_\mu|\}$.

\begin{lem}\label{lem-charge}
(i) If $c$ is such that $|c|\le1$, then $(C-c){\rm VIR}_K$ is a nontrivial ideal of ${\rm VIR}_K$, and the quotient ${\rm VIR}_K^c={\rm VIR}_K/(C-c){\rm VIR}_K$ is an admissible vertex $K$-algebra, whose field-state homomorphism $fs:{\rm Fld}({\rm VIR}_K^c) \to{\rm VIR}_K^c$ is a bijection;

(ii) If $c$ is invertible in $K$ and $|c|>1$, then $(C-c){\rm VIR}_K={\rm VIR}_K$;
\end{lem}

\begin{proof}
(i) Consider the Banach subspace of ${\rm VIR}_K$
 $$
W=\{f=\sum_\mu \lambda_\mu L^\nu \big\vert \lambda_\mu\in K\ {\rm and\ } ||\lambda_\mu||\to0\ {\rm as\ } \langle\mu\rangle\to\infty\}\ , 
 $$
where the sum is taken over sequences of non-negative integers $\mu=(\mu_1,\mu_2,\ldots)$ as above, and $||f||=\max_\mu\{||\lambda_\mu||\}$. Then there is a canonical isometric isomorphism ${\rm VIR}_K\toisom W\{C\}$, and we can work with $W\{C\}$ instead of ${\rm VIR}_K$. The multiplication by $C-c$ on $W\{C\}$ is the composition of the isometric automorphism $W\{C\} \to W\{C\}$ that is identical on $W$ and takes $C$ to $C-c$ and of the multiplication by $C$. This implies the claim.

(ii) If $c$ is invertible and $|c|>1$, the series $1+\frac{C}{c}+\frac{C^2}{c^2}+\ldots$ lies in $K\{C\}$ and its product with $1-\frac{C}{c}$ is equal to one. Since $K\{C\} \subset {\rm VIR}_K$, the claim follows.
\end{proof}

The algebra ${\rm VIR}_K^c$ is called the {\it universal Virasoro vertex $K$-algebra with central charge $c$}.

 \medskip
 
\centerline{\it Universal affine vertex $K$-algebra}
 \medskip
 
Let $\gothg$ be a Lie $R$-algebra, free of finite rank over $R=\bfZ[\frac{1}{N}]$ for an integer $N\ge1$, and assume $\gothg$ is provided with a non-degenerate symmetric $R$-bilinear form $\gothg\times\gothg\to R: (a,b)\mapsto (a|b)$. The latter gives rise to the associated Kac-Moody affinization $\widehat\gothg=\gothg[t^{\pm1}] + R{\rm K}$. defined by the relations
 $$
[at^m,bt^n] = [a,b] t^{m+n} + m\delta_{m,-n} (a|b){\rm K},\ [{\rm K},at^m]=0
 $$
for all $a,b\in\gothg$ and $m,n\in\bfZ$. For each $a\in\gothg$, consider the formal distribution $a(z)= \sum_{n\in\bfZ} (a t^n) z^{-n-1}$. Then the above relations can be written in the equivalent form
 $$
[a(z),b(w)]=[a,b](w)\delta(z-w) + (a|b) \partial_w \delta(z-w){\rm K},\ [{\rm K},a(z)]=0\ . 
 $$
It follows that $(z-w)^2 [a(z),b(w)]=0$ and, therefore, $\calF=\{a(z)\}_{a\in\gothg} \cup\{{\rm K}\}$ is a collection of mutually local formal distributions. The property (L.1) for $\calF$ clearly holds. Since $-\partial_t(a(z))=\partial_z a(z)$, the property (L.2) holds for $T=-\partial_t$. The above equality implies that $a(z)_{(0)}b(z)=[a,b](z)$, $a(z)_{(1)}b(z)=(a|b){\rm K}$, and $a(z)_{(i)}b(z)=0$ for $i\ge2$. This implies that the property (L.3) also holds. Thus, by Corollary \ref{cor-Lievertex}, we get a vertex $R$-algebra $V(\widehat\gothg,\calF)$ called the {\it universal affine vertex algebra} associated to the pair $(\gothg,(.|.))$. Notice that the field-state homomorphism ${\rm Fld}(V(\widehat\gothg,\calF)) \to V(\widehat\gothg,\calF)$ is a bijection.

The annihilation subalgebra of $\widehat\gothg$ is $\widehat\gothg_{-} = \gothg[t]$. Fix a basis $\{e_1,\ldots,e_p\}$ of $\gothg$ over $R$ and, for each $1\le i\le p$ and $n\ge1$, set $x_{i,n}=e_i t^{-n}\in\widehat\gothg$.  Then $V(\widehat\gothg,\calF)$ is a free $R$-module and, as an $R$-module, it is isomorphic to the ring of polynomials $R[{\rm K},x_{i,n}]_{1\le i\le p,n\ge1}$. A basis of the latter over $R$ is formed by the monomials ${\rm K}^m\cdot x^\mu$, where $m\ge0$,  $\mu=(\mu_1,\ldots,\mu_p)$, each $\mu_i$ is an infinite sequence of non-negative integers $(\mu_{i,1},\mu_{i,2},\dots)$ all but finitely many of which are equal to zero, and $x^\mu=x^{\mu_1}_1\cdot\ldots \cdot x^{\mu_p}_p$ with $x^{\mu_i}_i=\prod_{n=1}^\infty x_{i,n}^{\mu_{i,n}}$.

Assume now that $N$ is invertible in the commutative Banach ring $K$. By Lemma \ref{lem-fsurj}(i), $V_K (\widehat\gothg,\calF)= V(\widehat\gothg,\calF) \hatotimes_R K$ is an admissible vertex $K$-algebra, called the {\it universal affine vertex $K$-algebra associated to} $(\gothg,(\ |\ ))$. One has
 $$
V_K (\widehat\gothg,\calF) = \{f=\sum \lambda_{m,\mu} {\rm K}^m x^\mu \big\vert \lambda_{m,\mu}\in K\ {\rm and\ } \lambda_{m,\mu}\to0\ {\rm as\ } m+\langle\mu\rangle\to\infty\}\ , 
 $$
where $m$ and $\mu$ are as above and $\langle\mu\rangle = \sum_{i=1}^p \sum_{n=1}^\infty in\mu_{i,n}$. The following statement is verified in the same way as Lemma \ref{lem-charge}.

\begin{lem}\label{lem-chargeaff}
(i) If $|k|\le1$, then $({\rm K}-k)V_K (\widehat\gothg,\calF)$ is a nontrivial closed subspace of $V_K (\widehat\gothg,\calF)$, and the quotient $V^k_K (\widehat\gothg)=V_K (\widehat\gothg,\calF)/({\rm K}-k)V_K (\widehat\gothg,\calF)$ is an admissible vertex $K$-algebra whose field-state homomorphism $fs:{\rm Fld}(V^k_K (\widehat\gothg)) \to V^k_K (\widehat\gothg)$ is a bijection;
		
(ii) if $k$ is invertible in $K$ and $|k|>1$, then $({\rm K}-k)V_K (\widehat\gothg,\calF)=V_K (\widehat\gothg,\calF)$.
		\hfill$\square$
\end{lem}

The algebra $V^k_K (\widehat\gothg)$ is called the {\it universal affine vertex $K$-algebra of level $k$}.

For example, assume that $K=R=\bfZ$ and $\gothg$ is the commutative Lie algebra $\bfZ$, provided with the bilinear form $(a|b)=ab$. Then $V^1_K(\widehat\gothg)$ is canonically isomorphic to the free boson vertex $K$-algebra $B_K$.

In all of the above examples the Lie algebra in question is naturtally
$\bfZ$-graded. Namely, the Virasoro algebra (resp. affine Lie algebra) is graded by letting $\deg L_n =n$ (resp. $\deg at^n =n$ for $a\in \gothg$). These gradings induce $\bfZ$-gradings on the corresponding universal vertex $K$-algebras.
Since the highest component in these gradings of a singular vector is
a singular vector, it follows that the quotient of these vertex $K$-algebras by the maximal graded ideal is a simple vertex $K$-algebra. The same argument shows that the free boson and fermion vertex $K$-algebras are simple. 

\subsection{Vertex $K$-algebras and Lie conformal $K$-algebras}\label{subsec-lieconformal}

For a Banach $K$-module and a positive real number $r$, one denotes by $V \{r^{-1}\lambda\}$ the Banach $K$-module of formal power series $v=\sum_{n=0} \lambda^n v_n$ with $v_n\in V$ and $||a_n|| r^n\to0$ as $n\to\infty$, provided with the Banach norm $||v||=\max_n ||v_n||$. Notice that there is an isometric isomorphism $K\{r^{-1}\lambda\} \hatotimes_K V \toisom V \{r^{-1}\lambda\}$. If $K$ is a non-Archimedean field, $K\{r^{-1}\lambda\}$ is the algebra of functions analytic on the closed disc in the affine line over $K$ of radius $r$ with center at zero. 

\begin{defin}\label{defn-lieconformal}
A {\it Lie conformal $K$-algebra of radius $r$} is a Banach $K$-module $L$, provided with a bounded $K$-linear endomorphism $T:L\to L$ and a bounded $K$-linear homomorphism ({\it $\lambda$-bracket})
 $$
L\hatotimes_K L \to L\{r^{-1}\lambda\}: a\otimes b \mapsto [a_\lambda b]
 $$
with the following properties:
\begin{itemize}
	\item [(L.1)]
({\it sesquilinearity}) $[Ta_\lambda b] = -\lambda[a_\lambda b]$, $[a_\lambda Tb] = (\lambda+T) a_\lambda b]$;
    \item[(L.2)]
({\it skew-symmetry}) $[b_\lambda a] = - [a_{-\lambda-T} b]$;
    \item[(L.3)]
({\it Jacobi identity}) $[a_\lambda[b_\mu c]] = [[a_\lambda b]_{\lambda+\mu} c] + [b_\mu[a_\lambda c]]$.
\end{itemize}
\end{defin}

The right hand side in (L.2) should be understood as follows: if $[a_\lambda b] = \sum_{n=0}^\infty \lambda^n x_n$, then $[a_{-\lambda-T} b] = \sum_{n=0}^\infty (-\lambda-T)^n x_n$. 

Suppose now that our commutative Banach ring $K$ contains the field of rational numbers $\bfQ$. Then the valuation on $\bfQ$, induced by that on $K$, is either $p$-adic for a prime number $p$, or trivial. In the former case, we set $r_p=|p|^{\frac{1}{p-1}}$, and in the latter case, we set $p=1$ and $r_1=1$. In the following statements, $V$ is a vertex $K$-algebra.

\begin{prop}\label{prop-Lconf}
The space of fields ${\rm Fld}(V)$ is a Lie conformal $K$-algebra of radius $r_p$ with respect to the derivation $\partial_z$ and the $\lambda$-bracket, defined by
 $$
[\varphi(w)_\lambda\psi(w)] = {\rm Res}_z(e^{\lambda(z-w)} [\varphi(z),\psi(w)]) = \sum_{n=0}^\infty \frac{\lambda^n}{n!} (\varphi(w)_{(n)} \psi(w)) \ . 
 $$
\end{prop}

\begin{proof}
Step 1. {\it The second equality is true.} Indeed, the right hand side is equal to
$
\sum_{n=0}^\infty \frac{\lambda^n}{n!} {\rm Res}_z((z-w)^n [\varphi(z),\psi(w)])
$
and, by the definition in Lemma \ref{lem-localrep}(i), the latter is equal to $\sum_{n=0}^\infty \frac{\lambda^n}{n!} (\varphi(z)_{(n)} \psi(z))$.
\medskip

Step 2. {\it The sum on the right hand side of the above equality lies in ${\rm Fld}(V)\{r_p^{-1} \lambda\}$.} Indeed, if $p=1$, then $|n!|=1$, and since the fields $\varphi(z)$ and $\psi(z)$ are mutually local, then $||\varphi(z)_{(n)} \psi(z)||\to0$ as $n\to\infty$. This implies that the coefficient at $\lambda^n$ tends to zero as $n$ goes to infinity. Suppose $p$ is a prime. Then the degree of $p$ in the integer $n!$ is at most $[\frac{n}{p}] + [\frac{n}{p^2}] + \ldots \le \frac{n}{p-1}$. It follows that $|n!| \ge |p|^{\frac{n}{p-1}}$ and, therefore, $\frac{r_p^n}{|n!|} \le 1$. Since the fields $\varphi(z)$ and $\psi(z)$ are mutually local, the claim follows.
 \medskip
 
Step 3. {\it The equalities (L.1) hold.} Indeed, one has
 $$
[\partial_w  \varphi(w)_\lambda\psi(w)] = e^{-\lambda w} {\rm Res}_z([e^{\lambda z} \partial_z \varphi(z),\psi(w)])\ .
 $$
Since $e^{\lambda z} \partial_z \varphi(z) = \partial_z(e^{\lambda z} \varphi(z)) - \lambda e^{\lambda z} \varphi(z)$ and ${\rm Res}_z (\partial_z(e^{\lambda z} \varphi(z))) = 0$, the first equality follows. Furthermore, one has
 $$
[\varphi(w)_\lambda \partial_w \psi(w)] = {\rm Res}_z (e^{\lambda z} [\varphi(z), e^{-\lambda w} \partial_w \psi(w)])\ . 
 $$
Since $e^{-\lambda w} \partial_w \psi(w) = \lambda e^{-\lambda w} \psi(w) + \partial_w (e^{-\lambda w} \psi(w))$, the second equality follows.
 \medskip
 
Step 4. {\it The equality (L.2) holds.} Indeed, $[\varphi(z)_{-\lambda-\partial_z} \psi(z)]$ is equal to
\begin{eqnarray}
&\ & \sum_{n=0}^\infty \frac{(-\lambda-\partial_z)^n}{n!} (\varphi(z)_{(n)} \psi(z))  = \sum_{n=0}^\infty (-1)^n \sum_{i=0}^n \frac{1}{i!} \lambda^i \partial_z^{(n-i)} (\varphi(z)_{(n)} \psi(z)) {}\nonumber\\
&=& \sum_{i=0}^\infty \frac{\lambda^i}{i!} \sum_{m=0}^\infty (-1)^{m+i} \partial^{(m)}_z (\varphi(z)_{(m+i)} \psi(z))\ .
{}\nonumber
\end{eqnarray}
By Corollary \ref{cor-skew}(ii), the inner sum is equal to $-\psi(z)_{(i)} \varphi(z)$ and, therefore, the whole expression is equal to $-[\psi(z)_\lambda\varphi(z)]$.
 \medskip
 
Step 5. {\it The equality (L.3) holds.} Indeed, by Step 1, the left hand side of (L.3) is equal to
 $$
[\varphi(z)_\lambda[\psi(z)_\mu\chi(z)]] = {\rm Res}_z {\rm Res}_y (e^{\lambda(z-w)+\mu(y-w)} [\varphi(z),[\psi(y),\chi(w)]])\ .
 $$
By the classical Jacobi identity, one has
 $$
[\varphi(z),[\psi(y),\chi(w)]] = [[\varphi(z),\psi(y)],\chi(w)] + [\psi(y),\varphi(z),\chi(w)]]\ . 
 $$
Thus, the left hand side of (l.3) is a sum of
 $$
{\rm Res}_y (e^{(\lambda+\mu)(y-w)} [{\rm Res}_z (e^{\lambda(z-y)} [\varphi(z),\psi(y)]),\chi(w)]) \ ,
 $$
which is equal to $[[\varphi(z)_\lambda\psi(z)]_{\lambda+\mu} \chi(z)]$, and of
 $$
{\rm Res}_y (e^{\mu(y-w)} [\psi(y), {\rm Res}_z (e^{\lambda(z-w)} [\varphi(z), \chi(w)])])\ , 
 $$
which is equal to $[\psi(z)_\mu [\varphi(z)_\lambda\chi(z)]]$. The required equality follows.
\end{proof}

\begin{cor}\label{cor-Lconf}
In $V$ is admissible, the subspace $V'$ is a Lie conformal $K$-algebra of radius $r_p$ with respect to the translation operator $T$ and the $\lambda$-bracket, defined by
 $$
[a_\lambda b] = \sum_{n=0}^\infty \frac{\lambda^n}{n!} (a_{(n)} b)\ . 
 $$
\end{cor}

\begin{proof}
Since $V$ is admissible, the space-field homomorphism $fs: {\rm Fld}(V)\to V'$ is an isomorphism of Banach $K$-modules. The statement follows from Proposition \ref{prop-Lconf} and Lemma \ref{lem-basic}.
\end{proof}

\begin{exam}
Let $V$ be a commutative vertex $K$-algebra. Then ${\rm Fld}(V)\subset V[[z]]$, and since $\varphi(z)_{(n)}\psi(z)=0$ for all $\varphi(z),\psi(z)\in V[[z]]$ and $n\ge0$, it follows that the $\lambda$-bracket on the associated Lie conformal $K$-algebra is zero.	
\end{exam}

\end{document}